\documentclass[10pt]{article}
\usepackage[T1]{fontenc}
\usepackage[utf8]{inputenc}
\usepackage{color}
\usepackage{courier}
\usepackage{amsmath}
\usepackage{amsthm}
\usepackage{amssymb}
\usepackage{mathrsfs}
\usepackage{bbm}
\usepackage{graphicx}
\usepackage{epstopdf}
\usepackage{subfigure}
\usepackage{enumerate}
\usepackage{enumitem}
\usepackage{fancyhdr}
\usepackage{tikz}
\usepackage{makeidx}
\usepackage{thmtools}
\usepackage{thm-restate}
\usepackage{hyperref}
\usepackage{cleveref}
\newtheorem{theorem}{Theorem}
\newtheorem{proposition}[theorem]{Proposition}
\newtheorem{lemma}[theorem]{Lemma}
\newtheorem{corollary}[theorem]{Corollary}
\newtheorem*{theoremA}{Theorem A}
\newtheorem*{theoremAp}{Theorem A'}
\newtheorem*{theoremB}{Theorem B}
\newtheorem*{theoremBp}{Theorem B'}
\newtheorem*{theoremC}{Theorem C}
\newtheorem*{theoremD}{Theorem D}

\newcommand{\Z} {{\rm Z}\!\!{\rm Z}}

\newcommand{\T}{\mathrm{T}}
\newcommand{\CDiam}{\mathcal{D}}
\newcommand{\FM}{\mathcal{FM}}
\newcommand{\diff}{\operatorname{Diff}}

\begin{document}

\begin{titlepage}

\vskip0.5truecm

\begin{center}

{\LARGE \bf Homotopically unbounded disks for generic surface diffeomorphisms}

\end{center}

\vskip  0.4truecm

\centerline {{\large Salvador Addas-Zanata and Andres Koropecki}}

\vskip 0.2truecm

\centerline { {\sl Instituto de Matem\'atica e Estat\'\i stica }}
\centerline {{\sl Universidade de S\~ao Paulo}}
\centerline {{\sl Rua do Mat\~ao 1010, Cidade Universit\'aria,}} 
\centerline {{\sl 05508-090 S\~ao Paulo, SP, Brazil}}

\centerline{{\sl and}}

\centerline { {\sl Instituto de Matem\'atica e Estat\'\i stica }}
\centerline {{\sl Universidade Federal Fluminense}}
\centerline {{\sl Rua Professor Marcos Waldemar de Freitas Reis s/n}} 
\centerline {{\sl 24210-201 Niteroi, RJ, Brazil}}
 
\vskip 0.4truecm

\begin{abstract}

In this paper we consider closed orientable surfaces $S$ of positive genus and $C^r$-diffeomorphisms $f:S\rightarrow S$ isotopic to the identity ($r\geq 1)$. The main objective is to study periodic open topological disks which are homotopically unbounded (i.e. which lift to unbounded connected sets in the universal covering). We show that these disks are not uncommon, and are related to important dynamical phenomena. We also study the dynamics on these disks under certain generic conditions.
Our first main result implies that for the torus (or for arbitrary surfaces, with an additional condition) if the rotation set of a map has nonempty interior and is not locally constant, then the map is $C^r$-accumulated by diffeomorphisms exhibiting periodic homotopically unbounded disks. Our second result shows that $C^r$-generically, if the rotation set has nonempty interior (plus an additional hypothesis if the genus of $S$ is greater than $1$) a maximal periodic disk which is unbounded and has a rational prime ends rotation number must be the basin of some compact attractor or repeller contained in the disk. As a byproduct we obtain results describing certain periodic components of the complement of the closure of  stable or unstable manifolds of a periodic orbit in the $C^r$-generic setting.

\end{abstract} 

\vskip 0.5truecm


\vskip 0.8truecm

\noindent{\bf e-mail:} sazanata@ime.usp.br and ak@id.uff.br

\vskip 0.7truecm


\vfill
\hrule
\noindent{\footnotesize{The authors were partially supported 
by CNPq, grants: 306348/2015-2 and 305612/2016-6. SAZ was partially supported by FAPESP and AK was partially supported by FAPERJ.}}

 \end{titlepage}



\section{Introduction}

Periodic open disks are important objects for homeomorphisms of closed
orientable surfaces. In the positive genus and area preserving case, under
some natural conditions, their existence implies non ergodicity with respect
to Lebesgue measure \cite{korotal}. Several results have been proved in this
direction, see \cite{invent}, \cite{korotal}, \cite{c1epsilon}, \cite{bruno}
and \cite{final}. In the $C^r$-generic area-preserving setting, with $r$ large enough, the existence of elliptic periodic points gives rise to periodic disks (thus called ``elliptic islands'') due to the KAM phenomenon \cite{moser-kam}.

When the surface is the torus, there is a relationship between the instability of rotation sets and the existence of such disks: Recall from \cite{misiu} that for a homeomorphism $f\colon \T^2\to \T^2$ in the isotopy class of the identity, the rotation set of a lift $\tilde f\colon \mathbb{R}^2\to \mathbb R^2$ of $f$ is the (compact, convex) set $\rho(\tilde f)$ consisting of all possible limits of the form
$$\lim_{k\to \infty} \frac{\tilde f^{n_k}(z_k)-z_k}{n_k}$$ 
where $z_k\in \mathbb{R}^2$ and $n_k\to \infty$. 
In \cite{final} it is proved that for a $C^\infty $%
-generic area-preserving one-parameter family $f_t$ of diffeomorphisms of the torus, with $t$ belonging to some interval of the form $%
[t^{*}-1,t^{*}+1],$ if the rotation set $\rho (\widetilde{f}_t)$ always has interior (where $\widetilde f_t$ is a lifted family) and there is a rational vector $\rho\in \mathbb{Q}^2$ in the boundary of $\rho(f_{t^*})$ such that $\rho$ is in the interior of $\rho(f_t)$ for arbitrarily small values of $t-t^{*}>0,$ then  $t^{*}$
is accumulated by parameters $t_i$ for which $f_{t_i}$ has generic elliptic
periodic points. So every time a rational vector enters the interior of the
rotation set of a family as above, the critical parameter $t^{*}$ is
accumulated by parameters where elliptic islands are present. In this setting, the elliptic islands are always ``homotopically bounded'', in the sense that their lifts to the universal covering are bounded (see below). This is relevant because, despite the usefulness of the rotation set in studying the dynamics of toral homeomorphisms, the dependence of the rotation set on the map, and the mechanisms by which this set varies are poorly understood.

The main motivations for the present paper are to show that for general (non-conservative) generic one-parameter families of surface diffeomorphisms, unbounded periodic open disks often play the role of elliptic islands, and to understand the dynamics on such unbounded disks.

To be more precise, let $S$ be a closed surface of positive genus, endowed with a Riemannian metric, and $\pi \colon \widetilde S\to S$ the universal covering map. We say that an open topological disk $D\in S$ is homotopically bounded (or unbounded) if some (hence any) connected component of $\pi^{-1}(D)$ is bounded (resp. unbounded) in the lifted metric.
Given $r\geq 0$, denote by $\diff_0^r(S)$ the space of $C^r$ diffeomorphisms from $S$ to $S$ which are isotopic to the identity, endowed with the $C^r$ topology. Any such $f$ lifts to a $C^r$ diffeomorphism $\tilde f\colon \tilde S \to \tilde S$ which commutes with all deck transformations (unique if the genus of $S$ is greater than $1$) called \emph{natural lift of $f$}. 

Our first result is somewhat analogous to the one stated for area-preserving maps of the torus, but in the non-conservative setting (with homotopically unbounded disks instead of elliptic islands). Recall that having a rotation set with nonempty interior is a $C^0$-open condition.

\begin{theoremA}
	Given $r\geq 2$, let $(f_t)_{t\in I}$ be a $C^r$-generic one-parameter family in $\diff^r_0(\mathrm{T}^2)$ with a lift $(\tilde{f}_t)_{t\in I}$ such that the rotation set $\rho(\tilde f_t)$ has nonempty interior for each $t\in I$. Suppose that $t_0<t_1$ are such that a rational vector $\rho\in \mathbb{Q}^2$ lies outside of $\rho(\tilde f_{t_0})$ but in the interior of $\rho(\tilde f_{t_1})$. Then there exists an open interval $I^\prime \subset [t_0,t_1]$ such that for all $t\in I^\prime $, the map $f_t$ has a periodic homotopically unbounded open topological disk with rotation vector $\rho $ which is the basin of an attracting or repelling periodic point. 
\end{theoremA}

An immediate consequence is that if the rotation set (as a subset of $\mathbb{R}^2$ modulo integer translations) has nonempty interior and is not locally constant at $f\in \diff^r_0(\T^2)$, then $f$ is $C^r$-accumulated by diffeomorphisms exhibiting homotopically unbounded periodic disks (which are basins).

A similar result can be stated in a surface $S$ of higher genus. In order to be precise, let us introduce a definition: we say that a periodic point $Q$ of $f\in \diff^r_0(S)$ has a full mesh if $Q$ is a hyperbolic periodic saddle with a stable branch $\alpha^s$ (i.e. a connected component of $W^s(Q)\setminus\{Q\}$) and an unstable branch $\alpha^u$ such that every lift of $\alpha^s$ has a topologically transverse intersection with every lift of $\alpha^u$ (see Section \ref{sec:mesh} for more details). When $S=\T^2$ and $r\geq 2$, it is known from \cite{c1epsilon} that the rotation set of a lift of $f$ has nonempty interior if and only if there exists a periodic point with a full mesh.

Let us denote by $\FM_0^r(S)\subset \diff^r_0(S)$ the space of diffeomorphisms for which there exists a contractible hyperbolic periodic saddle $Q$ with a full mesh. By contractible we mean that it lifts to a periodic point of some natural lift of $f$ (see Section \ref{sec:prelim}). It is not difficult to see that this set is open in $\diff^r_0(S)$. Moreover, if $r\geq 2$, we know from \cite{c1epsilon} that $\FM_0^r(\T^2)$ coincides with the diffeomorphisms which have a lift for which the origin belongs to the interior of the rotation set, and for surfaces of higher genus a topological characterization (in terms of the existence of certain types of periodic orbits) is given in \cite{bruno} (see Section \ref{sec:mesh}).


%

The rotation set of $f\in \diff^0_0(S)$ can be defined when $S$ is a surface of higher genus $g>1$, and in this case it is a subset of the first homology group $ H_1(S, \mathbb{R})\simeq \mathbb{R}^{2g}$, which depends only on $f$ and not on the choice of a particular isotopy from the identity to $f$ (see Section \ref{sec:prelim}). The same conclusion of theorem A holds in this setting for generic one-parameter families in the space $\FM^r_0(S)$: 

\begin{theoremAp}
	Suppose the genus of $S$ is $g>1$. Given $r\geq 2$, let $(f_t)_{t\in I}$ be a $C^r$-generic one-parameter family in $\FM^r_0(S)$. Suppose that $t_0<t_1$ are such that a rational vector $\rho\in H_1(S, \mathbb{Q}) \simeq \mathbb{Q}^{2g}$ lies outside of $\rho(f_{t_0})$ but in the interior of $\rho(f_{t_1})$. Then there exists an open interval $I^\prime\subset [t_0,t_1]$ such that for all $t\in I^\prime $, the map $f_t$ has a periodic homotopically unbounded open topological disk with rotation vector $\rho$, which is the basin of an attracting or repelling periodic point.
\end{theoremAp}

Recall that a topological attractor of $f$ is a compact set $K$ which has a neighborhood $U$ such that $f(\overline U)\subset U$ and $K=\bigcap_{n\geq 0} f^n(U)$. The basin (of attraction) of $K$ is the set $\bigcup_{n\in \Z} f^n(U)$. A topological repeller is defined similarly using $f^{-1}$ instead of $f$. 
Finally, a periodic open topological disk is maximally periodic if it is not properly contained in any other periodic open topological disk. Note that every homotopically unbounded periodic disk is contained in a maximally periodic homotopically unbounded disk. 

Our second main result describes the dynamics of (maximal) homotopically unbounded periodic disks when the rotation set has nonempty interior.  

\begin{theoremB} Given $r\geq 1$, let $f$ be a $C^r$-generic element of $\diff^r_0(\T^2)$, such that  the rotation set of a lift of $f$ has nonempty interior. If $D = f^n(D)$ is a homotopically unbounded maximally periodic open topological disk and its prime ends rotation number\footnote{see Section \ref{sec:primends}} is rational, then $D$ is a basin (of attraction or repulsion) for $f^n$.
\end{theoremB}

Thus, the ``unbounded'' part of $D$ actually converges under forward or backward iteration towards a compact
connected (hence bounded) attracting or repelling set. 

The topology of connected components of the lift to the universal covering of
homotopically unbounded maximally periodic open disks can be very
complicated. In the torus case, for instance, assume $D\subset {\rm T^2}$ is
such a disk. Then, $f^n(D)=D$ for some integer $n\geq 1.$ Fixed some lift of 
$f$ to the plane, denoted $\widetilde{f},$ there exists an integer pair $%
(a,b)$ such that for any connected component $\widetilde{D}$ of $\pi
^{-1}(D),$ $\widetilde{f}^n(\widetilde{D})=\widetilde{D}+(a,b).$ From
proposition \ref{geral}, for any periodic point $z\in {\rm T^2}$
with  a full mesh, either $D\cap W^s(z)\neq \emptyset ,$ and in this
case $\partial D=\overline{W^u(z)}$
from the maximallity of $D,$ or $D\cap W^u(z)\neq \emptyset $ and 
$\partial D=\overline{W^s(z)}.$ So, if $(a/n,b/n),$ the rotation vector
of $D,$ belongs to the interior of $\rho (\widetilde{f}),$ we could choose
the point $z$ which has a full mesh with rotation vector equal to $%
(a/n,b/n),$ see theorem \ref{casogeral}. And considering the lift of $f^n$
to the plane, given by $\widetilde{h}=\widetilde{f}^n(\bullet )-(a,b),$ we
get that both $\widetilde{D}$ and any $\widetilde{z}\in \pi ^{-1}(z)$ are $
\widetilde{h}$-periodic. Thus, as $\widetilde{D}$ either intersects $W^s(
\widetilde{z})$ or $W^u(\widetilde{z}),$ the $\lambda $-lemma 
\cite{palis} implies that 
$\partial \widetilde{D}$ is equal to either $\overline{W^s(\widetilde{z})}$
or $\overline{W^u(\widetilde{z})},$ which are closed connected equivariant
subsets (invariant under integer translations) of the plane. In particular, $
\widetilde{D}$ is an open disk which intersects all fundamental domains of
the torus and clearly, as $D=\pi (\widetilde{D})$ is a disk, $\widetilde{D}%
\cap (\widetilde{D}+(t,s))=\emptyset $ for all integer pairs $(t,s)\neq
(0,0).$ So, all connected components of $\pi^{-1}(D)$ have the same 
common boundary!   

The other possibility is when the rotation vector of $D$ belongs to the
boundary of the rotation set. In this case, from the bounded displacement
condition satisfied with respect to vectors in the boundary of the rotation
set (see \cite{davalos}, \cite{london} and \cite{fabio}), $D$ can not be as
above. It has to be bounded in at least some direction.

 About the hypotheses of the above theorem (and of theorem B'), we do not know
 whether the condition on the prime ends rotation number is necessary: 
\newtheorem*{question*}{Question}
\begin{question*}
	Is there a homeomorphism $f$ of a closed surface with an $f$-invariant homotopically unbounded open topological disk $D$ having an irrational prime ends rotation number?
\end{question*}
Even if we assume that $f$ is smooth, the surface is $\T^2$, and the rotation set of $f$ has nonempty interior, we do not know the answer to this question.

A version of theorem B holds in higher genus surfaces if one considers elements of $\FM^r_0(S)$:
\begin{theoremBp} Let $S$ be a closed orientable surface of positive genus, and given $r\geq 1$ let $f$ be a $C^r$-generic element of $\FM^r_0(S)$. Then any  homotopically unbounded maximally periodic open topological disk $D=f^n(D)$ with a rational prime ends rotation number is a basin (of attraction or repulsion) for $f^n$.
\end{theoremBp}

Theorems B and B' are a consequence of the following more general fact:
\begin{theoremC}
	\label{atra} For any $r\geq 1$, if $f\in \diff^r_0(S)$ is a 
$C^r$-generic diffeomorphism of a closed surface $S$ of positive genus, $z$ is a contractible periodic point with a full mesh, $D = f^n(D)$ is a periodic connected component of the complement of $\overline{W^s(z)}$ and the prime ends rotation number of $f^n$ in $D$ is rational, then $D$ is a basin of attraction. 
\end{theoremC}
Of course, the analogous result holds replacing $W^s(z)$ with $W^u(z)$, in which case $D$ is a basin of repulsion. In addition, under the hypotheses of theorem C we are able to give a rather detailed description of the prime ends dynamics: there are finitely many periodic prime ends, some of which are accessible hyperbolic saddles and the remaining are sources. For the saddles, one branch of the unstable manifold is contained in $D$ and the whole stable manifold belongs to $\partial D$ and consists of accessible points. See theorem \ref{main1}.

Note that in theorem C above it is irrelevant whether $D$ is homotopically unbounded or not.
To prove theorems B and B' from this result, one shows that a maximally periodic homotopically unbounded open topological disk must always be a connected component of $S\setminus \overline{W^s(z)}$ or of $S\setminus  \overline{W^u(z)}$ for a contractible periodic point $z$ with a full mesh. 

A final observation about theorems B, B' and C is the following: In order to prove that diffeomorphisms of the torus isotopic to the identity whose rotation sets have non empty interior, have periodic points with a full mesh, theorem \ref{casogeral} needs $C^{1+\epsilon}$ differentiability for some $\epsilon>0$. The same is true for theorem \ref{casogeral2}, in the higher genus case, where for $C^{1+\epsilon}$ diffeomorphisms isotopic to the identity, the presence of a fully essential system of curves (see Definition \ref{thesystem}) implies existence of a contractible periodic point with a full mesh. But in the statements of theorems B, B' and C, the generic $C^1$ case is also considered. Nevertheless, it is easy to deal with it, because the subset of $C^2$ diffeomorphisms is dense in the set of $C^1$ diffeomorphisms. So, in the torus case, from theorem \ref{casogeral}, a $C^2$ diffeomorphism whose rotation set has non empty interior (remember that this is a $C^0$ open condition) will have a periodic point with a full mesh, as we said, a $C^1$ open condition. So, the $C^1$-genericity in theorem B is understood.

In the higher genus case, theorems B' and C already assume the existence of contractible periodic points with a full mesh, so there is nothing to add. But one could ask, if instead of assuming the existence of contractible points with full mesh, a $C^1$-generic $f$ with a fully essential system of curves were considered. What would happen? It turns out that the existence of a fully essential system of curves means that the map $f$ has a finite number of periodic orbits with certain special topological properties. In the $C^1$-generic case, these orbits are non-degenerate and thus persist under small perturbations (the above mentioned topological properties of these orbits persist as well, again see Definition \ref{thesystem}). Thus if $f^*$ is a $C^2$ diffeomorphism sufficiently close to $f$, then it also has a fully essential system of curves and so, theorem \ref{casogeral2} applies to $f^{*},$ assuring the existence of a contractible periodic point with a full mesh. As we said above, this is a $C^1$ open condition. So, if instead of assuming the existence of saddles with full mesh in theorems B' and C, we considered generic $C^r$ (for any $r\geq 1$) diffeomorphisms having fully essential systems of curves, the same conclusions would hold. For $r\geq 2$, this would be a trivial consequence of theorem \ref{casogeral2} and for $r=1$, it would follow from the above discussion.         

Using the same ideas from theorem C, one can show a general result in the contractible case:

\begin{theoremD} Let $M$ be a closed orientable surface and fix any $r\geq 1.$ For any $C^r$-generic diffeomorphism 
$f \in \diff^r_0(M),$ which we assume to have a periodic point with a full mesh if $M$ has positive genus, let 
$z$ be a hyperbolic $n$-periodic saddle point such that $\overline{W^s(z)}$ is contractible. If $D$ is a periodic connected
component of the complement of $ \cup_{i=0}^{n-1}f^i(\overline{W^s(z)})$
 which has rational prime ends
rotation number at one of its boundary components (therefore, at all of them), 
then $D$ is a basin of attraction ( if instead of stable manifold of $z$, we 
considered unstable manifold, then $D$ would be the basin of a repeller).
\end{theoremD}

A nice corollary of the above theorem is the following:

\vskip0.2truecm

{\bf Corollary of theorem D.}
{\it Fixed any $r\geq 1$ and any $C^r$-generic diffeomorphism $f$ of the sphere
without homoclinic points, for any hyperbolic periodic saddle $z$ of $f$, each
branch at $z$ does not accumulate, neither on itself nor on any of the other
three branches at $z$ (neither on the whole orbit of all the branches). Moreover, the unstable branches accumulate on connected
attractors and the stable ones, accumulate on connected repellers. 
}

\vskip0.2truecm

This result appears in \cite{crotrepuj} under other hypotheses for planar diffeomorphisms and it is a first step in showing that zero entropy mildly dissipative diffeomorphisms of the plane are either infinitely renormalizable or Morse-Smale.

This paper is organized as follows. In Section 2 we present some results we use and definitions. In Section 3, the results proved by Pixton for planar surfaces in \cite{pixton} are generalized to other surfaces with some additional hypotheses. In Section 4 we prove the main technical result of the paper, in Section 5 we extend some generic type results that hold in the space of diffeomorphisms to the parameter space of generic families, and in Section 6 proofs of theorems A,A',B,B',C and D appear.
	
\section{Preliminaries and additional results}
\label{sec:prelim}

\subsection{Notation and some definitions}

Let $S$ be a closed orientable surface. We denote by  $\pi :\widetilde{S}\rightarrow S$ the universal covering map of $S$, and by $Deck(\pi)$ the group of deck transformations of $\pi$. 

We assume that surfaces are endowed with a Riemannian metric, and their universal coverings with the lifted metrics. For any open topological disk $D\subset {S}$ we define $\CDiam(D) $ as the diameter of any connected component of $\pi^{-1}(D)$ (this does not depend on the choice of the component). We say $D$ is homotopically unbounded when $\CDiam(D)=\infty ,$ otherwise we say it is homotopically bounded.

As we already defined, for any $r\geq 0$, $\diff_0^r(S)$ is the set of $C^r$ diffeomorphisms of $S$ isotopic to the identity (clearly, when $r=0$, we mean homeomorphisms).

A compact subset $K$ of $S$ is said to be contractible if there exists a closed topological disk $\overline{D}\subset S$ such that $K\subset int(\overline{D}%
).$ 
We will make frequent use of the following observation: if $K$ is connected and
contractible, then the connected components
of $\pi^{-1}(K)$ are bounded and homeomorphic to $K.$

\subsection{Prime ends compactification of open disks}
\label{sec:primends}

If $D$ is an open topological disk of an oriented surface such that $\partial D$ is a Jordan curve and $f$ is an orientation preserving homeomorphism of that
surface which satisfies $f(D)=D,$ it is easy to see that $f:\partial
D\rightarrow \partial D$ is conjugate to a homeomorphism of the circle, and
so a real number $\rho (D)=rotation$ $number$ $of$ $f\mid _{\partial D}$ can
be associated to this problem. Clearly, if $\rho (D)$ is rational, then
there exists a periodic point in $\partial D$ and if it is not, then there
are no such points. This is known since Poincar\'e. The difficulties arise
when we do not assume $\partial D$ to be a Jordan curve.

The prime ends compactification is a way to attach to $D$ a circle called
the circle of prime ends of $D,$ obtaining a space $D\sqcup S^1$ with a
topology that makes it homeomorphic to the closed unit disk. If, as above,
we assume the existence of an orientation preserving homeomorphism $f$ such
that $f(D)=D,$ then $f\mid _D$ extends to $D\sqcup S^1.$ The prime ends
rotation number of $f$ in $D,$ still denoted $\rho (D),$ is the usual
rotation number of the orientation preserving homeomorphism induced on $S^1$
by the extension of $f\mid _D.$ But things may be quite different in this
setting. In full generality, it is not true that when $\rho (D)$ is
rational, there are periodic points in $\partial D$ and for some examples, $%
\rho (D)$ is irrational and $\partial D$ is not periodic point free. Here we
refer to \cite{mather} and \cite{koropatmey} for definitions of prime ends,
prime chains, end-cuts, cross-cuts, cross-sections, principal points, as
well as some important theorems. 

Let us state an addendum to theorem C describing the prime ends dynamics, which is part of our proof:

\begin{theorem}
	\label{main1} Under the hypotheses of theorem C, if $D$ is a
	connected component of $(\overline{W^s(z)})^c$ and $\widehat{D}$ is the prime ends compactification of $D,$ then $\widehat{D}$ is a closed disk in whose boundary there are only finitely many periodic prime ends.
	Necessarily, some are accessible hyperbolic saddles and the remaining are
	sources. For the saddles, one branch of the unstable manifold is contained
	in $D$ and the whole stable manifold belongs to $\partial D$ and is made of
	accessible points.
\end{theorem}

\subsection{Rotation sets}


\paragraph{The torus.}

Let $\mathrm{T}^2=\mathbb{R}^2/\mathbb{Z}^2$ be the flat torus and let 
$\pi :\mathbb{R}^2\longrightarrow {\mathrm T^2}$ be the associated covering map.
Coordinates are denoted as $\widetilde{x}\in \mathbb{R}^2$
and $x\in {\mathrm T^2.\ }$ Note that the deck transformations of $\pi$ in this case are the integer translations on the plane.

We denote by $\diff_0^r(\mathbb{R}^2)$ the set of lifts of elements from $\diff_0^r({\mathrm T^2})$ to the plane. Maps from 
$\diff_0^r(\mathrm{T^2})$ are denoted $f$ and their lifts to the plane are denoted $\widetilde{f}.$ Any such lift commutes with the deck transformations (i.e. with integer translations).

Given $f\in \diff_0^0({\mathrm T^2})$ and a lift $\widetilde{f}\in \diff_0^0(\mathbb{R}^2),$
 the Misiurewicz-Ziemian rotation set of $
\widetilde{f},$ $\rho (\widetilde{f}),$ can be defined as follows (see \cite
{misiu}): 
\begin{equation}
	\label{rotsetident}\rho (\widetilde{f})=\bigcap_{{\ 
			\begin{array}{c}
				i\geq 1 \\  
			\end{array}
	}}\overline{\bigcup_{{\ 
				\begin{array}{c}
					n\geq i \\  
				\end{array}
		}}\left\{ \frac{\widetilde{f}^n(\widetilde{z})-\widetilde{z}}n:\widetilde{z}
		\in \mathbb{R}^2 \right\}} 
\end{equation}

This set is a compact convex subset of $\mathbb{R}^2$ (see \cite
{misiu}), and it was proved in \cite{franksrat} and \cite{misiu} that all
points in its interior are realized by compact $f$-invariant subsets of $%
{\rm T^2,}$ which can be chosen as periodic orbits in the rational case. By
saying that some vector $\rho \in \rho (\widetilde{f})$ is realized by a
compact $f$-invariant set, we mean that there exists a compact $f$-invariant
subset $K\subset {\rm T^2}$ such that for all $z\in K$ and any $\widetilde{z}%
\in \pi^{-1}(z)$ 
\begin{equation}
	\label{deffrotvect}{\lim }_{n\rightarrow \infty }\frac{\widetilde{f}^n( 
		\widetilde{z})-\widetilde{z}}n=\rho . 
\end{equation}
Moreover, the above limit, whenever it exists, is called the rotation vector
of the point $z,$ denoted $\rho (f,z)$$.$

\paragraph{Surfaces of genus larger than one.}

When $S$ is a closed orientable surface of genus $g>1$, we may identify
the universal covering $\widetilde{S}$ with the Poincar\'e disk $\mathbb{D}$. As before, we denote by $\diff_0^r(\mathbb{D}%
)$ the set of lifts of elements from $\diff_0^r(S)$ to the Poincar\'e
disk. However here we only consider lifts which commute with all deck
transformations. It is well known that, given $f\in \diff_0^r(S),$ there
exists only one $\widetilde{f}:\mathbb{D}\to \mathbb{D}$ lifting $f$ and
having this property. It is called the natural lift of $f.$ One way to find
it, is to consider an isotopy $I$ from $Id:S\rightarrow S$ to $%
f:S\rightarrow S.$ The natural lift $\widetilde{f}$ is given by the endpoint
of the lift of the isotopy $I$ which starts at $Id:\mathbb{D}\to \mathbb{D}.$
As in the torus case, maps from $\diff_0^r(S)$ are denoted $f$ and their lifts
to $\mathbb{D}$ are denoted $\widetilde{f}.$

The notion of Misiurewicz-Ziemian rotation set introduced for the torus can be generalized to a surface $S$ of genus larger than one. In order to avoid unnecessary details, we refer to \cite{korotal} for a precise definition. Here we limit ourselves to stating that the rotation set of $f\in \diff_0^0(S)$ is a subset $\rho_{mz}(f)$ of the first homology group  $H_1(S,\mathbb{R})\simeq \mathbb{R}^{2g}$ (where $g$ is the genus of $S$), and that periodic orbits of $f$ have a corresponding rotation vector which is a rational element of the rotation set, as follows: given an isotopy from the identity to $f$ and a periodic point $p$ of least period $n$, the rotation vector of $p$ is $\alpha/n\in H_1(S, \mathbb{Q})$, where $\alpha\in H_1(S, \mathbb{Z})$ represents the homology class of the loop obtained by following the isotopy from $p$ to $f(p)$, then from $f(p)$ to $f^2(p)$, and so on until $f^n(p) = p$.

The set $\rho_{mz}(f)$ does not need to be convex for surfaces of higher genus, nor it is generally true that every rational element of its interior is realized by a periodic point. However, following \cite{bruno}, one has that these properties are true if $f\in \FM_0^r(S)$. Recall that this means that $f$ has a contractible hyperbolic periodic saddle (i.e. one that lifts to a periodic point of the natural lift of $f$) with a full mesh, and for smooth enough maps this is guaranteed by certain topological conditions, explained in the next Section.

\subsection{Existence of full mesh}
\label{sec:mesh}

Recall from the introduction that a hyperbolic periodic saddle point $Q$ of $f\in \diff_0^1(S)$ is said to have a full mesh if there exist branches $\alpha^u$ of $W^u(Q)$ and $\alpha^s$ of $W^s(Q)$ such that every lift of $\alpha^u$ to the universal covering $\tilde S$ has a topologically transverse intersection with every lift of $\alpha^s$. In other words, if $\tilde \alpha^u$ and $\tilde \alpha^s$ are any two lifts of $\alpha^u$ and $\alpha^s$, respectively, then $\tilde \alpha^u$ has a topologically transverse intersection with $g(\tilde \alpha^s)$ for every deck transformation $g$ of the universal covering. For instance, in the case $S=\T^2$, this means that $\tilde \alpha^s$ intersects every integer translation of $\tilde \alpha^u$.

For a precise definition of topologically transverse intersections, see definition 9, page 4 of \cite{c1epsilon}. In the present paper, however, we almost will not need it, because most of the time we consider generic maps, and for a $C^r$-generic diffeomorphism ($r\geq 1$) all intersections between stable and unstable manifolds of periodic points are $C^1$-transverse. We remark that the existence of a topologically transverse intersection between the stable and unstable manifolds of two hyperbolic saddles is a $C^r$-open condition for any $r\geq 1$.

It is not difficult to see that if $f\in \diff^1_0(\T^2)$ has a periodic point with a full mesh, then the rotation set of a lift of $f$ has nonempty interior.
In \cite{c1epsilon}, the following result was proved, which implies the converse of the previous claim in the $C^{1+\epsilon}$ case.

\begin{theorem}
\label{casogeral} Suppose $f$ belongs to $\diff_0^{1+\epsilon }({\rm T^2})$
for some $\epsilon >0$ and $(p/q,r/q)\in int(\rho (\widetilde{f})).$ Then,
for some integer $n\geq 1,$ $f$ has a $nq$-periodic hyperbolic saddle point 
$Q$ of rotation vector $(p/q,r/q)$ which has a full mesh.
\end{theorem}

For surfaces of higher genus, a generalization of the previous result appeared in \cite{bruno}. 
In order to state it, first we remember what an inessential subset of a closed orientable 
surface $S$ is: $K \subset S$ is inessential, if and only if, it is contained in an
homotopically bounded open topological disk of $S$.
 
And now, we introduce a definition from \cite{bruno}.

\begin{description}
	\item[Definition]  \label{thesystem} \ref{thesystem} {\bf Fully essential
		system of curves $\mathscr{C}:$} We say that $f\in \diff_0^0(S)$ is a
	homeomorphism with a fully essential system of curves $\mathscr{C}=\cup
	_{i=1}^k\gamma _i,$ if the following conditions are satisfied ($I$ is some fixed 
      isotopy from $Id$ to $f$):
	
	\begin{enumerate}
		\item  each $\gamma _i$ is a closed geodesic in $S$ and the complement of 
        $\cup _{i=1}^k\gamma _i$ only has inessential connected components $(k\geq 1)$;
		
		\item  for each $i\in \{1,\ldots ,k\},$ there is a $f$-periodic point $p_i$
		such that its trajectory under the isotopy $I$ is a closed curve freely
		homotopic to $\gamma _i$ with the correct orientation (the $\gamma
		_{i^{\prime }s}$ are oriented);
		
		\item  for every open intervals $B,F\subset \partial \mathbb{D}$, there
		exists an oriented simple arc $\widetilde{\alpha }\subset \pi ^{-1}(%
		\mathscr{C})$ formed by the concatenation of a finite number of oriented
		subarcs of extended lifts of geodesics in $\mathscr{C}$ and such that the
		initial point of $\widetilde{\alpha }$ is contained in $B$ and the final
		point belongs to $F$.
	\end{enumerate}
	
	{\bf Remark: }See \cite{bruno} for some comments on this definition. In
	particular, instead of the third condition, we could assume that for each $%
	\gamma _i,$ there are two $f$-periodic points $p_i^{+}$ and $p_i^{-}$ such
	that their trajectories under the isotopy are closed curves freely homotopic
	to $\gamma _i$ and to $-\gamma _i.$ It is not hard to see that this
	assumption implies the existence of a fully essential system of curves as
	stated above, see propositions 6 and 7 of \cite{bruno}.
\end{description}

\begin{theorem}
\label{casogeral2} Suppose $f$ belongs to $\diff_0^{1+\epsilon }(S)$ for some 
$\epsilon >0$ and it has a fully essential system of curves {\bf $\mathscr{C}%
.$} Then, $f$ has a contractible hyperbolic periodic saddle point $Q\in S$
with a full mesh.
\end{theorem}

The following discussion can be made assuming only topologically transverse
intersections. But as in most of this paper we are assuming generic
conditions, we will consider until the end of this subsection, that all
intersections between stable and unstable manifolds of periodic points are $%
C^1$-transverse. This simplifies several arguments.

In both cases, the torus and higher genus, it is possible to find finitely
many oriented closed curves in the surface, such that each curve is given by
the union of an arc in $W^u(Q),$ starting at $Q$ and ending at a point $z$
in $W^u(Q)\cap W^s(Q),$ and another arc contained in $W^s(Q),$
starting at $z$ and ending at $Q,$ with the property that the union of these
closed curves generates the homotopy group of the surface as a semi-group
(the orientation of each curve is given by the direction of the arrows in
the unstable and stable manifolds).

In particular, either for the torus or for higher genus surfaces, we get the
following:

\begin{proposition}
\label{bounddiam} Suppose $M$ is a closed
surface of positive genus and $f\in \diff_0^1(M)$ has a contractible
hyperbolic periodic saddle point with a full mesh. Then there is another
(possibly the same) hyperbolic periodic saddle $Q$ and compact arcs $\lambda
_u\subset W^u(Q)$ and $\lambda _s\subset W^s(Q)$ such that $Q\in \lambda
_s\cap \lambda _u$ and $M\setminus (\lambda _s\cup \lambda _u)$ is a union
of open topological disks whose covering diameters are uniformly bounded
from above by some constant $Max(f)>0$. Furthermore, $Q$
can be chosen with infinitely many different rotation vectors. 
In particular, given a bounded continuum $\widetilde{K}\subset \widetilde{M}$
and an integer $n\geq 1$ such that, 
$\widetilde{f}^n(\widetilde{K})=g(\widetilde{K})$ for some deck
transformation $g,$ then the $f^n$-periodic continuum $K\subset M$ given by $%
\pi (\widetilde{K})$ is actually contractible.
\end{proposition}

{\it Proof: }First, we need to prove the assertion about infinitely many
different rotation vectors for surfaces of genus larger than 1. For the
torus, we know that for each rational vector in the interior of the rotation
set, there exists a saddle with a full mesh with that rotation vector.

But in the higher genus context, it is unknown and maybe not even true that
for all rational vectors in the interior of the rotation set, there exists a
saddle with a full mesh with that rotation vector. Nevertheless, we can
consider the following construction (for more details, see lemma 18 and the 
proof of theorem 1 of \cite{bruno}): 
Let $Q_0$ be the contractible saddle which has a full
mesh. Picking some deck transformation $h$ whose homology class $[h]$ is non
trivial, we can build a rotational horseshoe at $Q_0$ so that symbol $"0"$
corresponds to the vertical rectangle that contains $Q_0$ and symbol $\;"1"$
corresponds to the vertical rectangle that moves in the direction of $h$ in
the universal covering, under some adequate iterate of $f.$ So, it is
straightforward that we can find periodic saddles $z_\alpha $ in this
horseshoe, such that varying $\alpha $ in an infinite set, the $z_{\alpha
^{\prime }s}$ have infinitely many different rotation vectors (infinitely
many rational multiples of $[h]$). And moreover, the stable (unstable)
manifold of $Q_0$ intersects the unstable (resp. stable) manifold of $%
z_\alpha $ (transversely). So, $\overline{W^u(Q_0)}=\overline{W^u(z_\alpha )}
$ and $\overline{W^s(Q_0)}=\overline{W^s(z_\alpha )},$ which imply that for
all $\alpha ^{\prime }s,$ there are compact arcs $\lambda _u^\alpha \subset
W^u(z_\alpha )$ and $\lambda _s^\alpha \subset W^s(z_\alpha )$ such that $%
z_\alpha \in \lambda _s^\alpha \cap \lambda _u^\alpha $ and $M\setminus
(\lambda _s^\alpha \cup \lambda _u^\alpha )$ is a union of open topological
disks.

Now, it is easy to prove the second part of the proposition: Assume $
\widetilde{K}$ is a continuum as in the statement. So, there exist $n>0$ and a 
deck transformation $g$ such that $\tilde{f}^n(\widetilde{K})=g(\widetilde{K})$.
This implies that all points in $K$ have the same well-defined rotation vector 
$[g]/n$ for $f$ (and $-[g]/n$ for $f^{-1}$). So, if we pick the saddle $Q$ equal
to some $z_{\alpha ^{*}}$ so that its rotation vector is different from the
rotation vector of $K,$ then as $\left( \lambda _u^{\alpha ^{*}}\cup \lambda
_s^{\alpha ^{*}}\right) $ cannot intersect $K,$ the proof is over. $\Box $
\vskip0.2truecm

\section{Some extensions of Pixton's theorems}

Throughout this section, $M$ denotes a surface, and some integer $r\geq 1$ is fixed.
We say that the surface $M$ is \emph{planar} if each connected component of $M$ is homeomorphic to a subset of the sphere. If $M$ is any surface, we say that a compact subset $K\subset M$ is planar if $K$ has a neighborhood in $M$ which is a planar surface. If $f\in \diff^r(M)$ is given, the orbit of a subset $Z\subset M$ is denoted $\mathcal{O}_f(Z)=\{f^n(z):z\in Z,\, n\in \Z\}$.

If $p$ is a hyperbolic periodic saddle point of $f\in \diff^r(M)$, a \emph{stable (or unstable) branch of $f$ at $p$} is a connected component of $W^s(p)\setminus \{p\}$ (resp. $W^u(p)\setminus\{p\}$). A periodic branch of $f$ is a stable or unstable branch of some hyperbolic periodic saddle point of $f$.

Also as a consequence of the local linearization of hyperbolic periodic saddles, one may find an arbitrarily small neighborhood $U$ of $p$ such that the local stable and unstable manifolds of $p$ in $U$ split $U$ into four ``local quadrants''. By the local stable (or unstable) manifold of $f$ in $U$ we mean the connected component $W^s_{loc}(p)$ (resp. $W^u_{loc}(p)$) of $W^s(p)\cap U$ (resp. $W^u(p)\cap U$) containing $p$.
We say that a set $Z$ accumulates on $p$ through a local quadrant $Q$ if $Z\cap Q$ accumulates on $p$. 
The local stable and unstable branches at $p$ adjacent to $Q$ are the connected components $Y_{loc}$ of $W^s_{loc}(p)\setminus\{p\}$ and $X_{loc}$ of $W^u_{loc}(p)\setminus\{p\}$ which lie in the boundary of $Q$ in $U$;  the stable and unstable branches adjacent to $Q$ are the sets $Y = \bigcup_{n\geq 0} f^{-kn}(Y_{loc})$ and $X=\bigcup_{n\geq 0} f^{kn}(X_{loc})$, where $k>0$ is an integer such that $f^k(Y_{loc}) \subset Y_{loc}$ and $f^{-k}(X_{loc}) \subset X_{loc}$.

If $Y$ is a stable branch and $X$ is an unstable branch of $p$, we say that \emph{a set $Z$ accumulates on a point $x\in X$ locally from the $Y$-side} if the following holds: given any compact arc $\gamma^u$ of $W^u(p)$ containing both $p$ and $x$ in its relative interior, and any  tubular neighborhood $U$ of $\gamma^u$, the connected component $U_Y$ of $U\setminus \gamma^u$  containing the initial local stable branch $Y_{loc}$ corresponding to $Y$ is such that $U_Y\cap Z$ accumulates on $x$. 

For future reference, we state the following fact, which is a straightforward consequence of the linearization of hyperbolic periodic saddles and the existence of horseshoes:

\begin{proposition}\label{Rectangles} Let $Y$ and $X$ be a stable and an unstable branch (respectively) of a hyperbolic periodic saddle of $f\in\mathrm{Diff}^r(M)$, and suppose $X$ meets $Y$ transversely. If $Z$ is a periodic arcwise connected set, then $Z$ intersects $X\cup Y$ provided that one of the following holds:
	\begin{itemize}
		\item[(1)] One of the branches $X$ or $Y$ is contained in $\overline{Z}$, or
		\item[(2)] $Z$ accumulates on $X$ from the $Y$ side, or
		\item[(3)] $Z$ accumulates on $Y$ from the $X$ side.
	\end{itemize}
	In particular, if $Z$ is an unstable (or stable) branch of a periodic saddle and one of the above cases holds, then $Z$ must intersect $Y$ (resp. $X$).
\end{proposition}

\begin{proposition}
	\label{contemramo}  Suppose $f\in \diff^r(M)$ and $K\subset M$ is an $f$-periodic continuum not reduced to a point. 
	
	\begin{enumerate}
		\item  If $p\in K$ is a hyperbolic periodic saddle, then at least one of the four branches at $p$ is contained in $K;$
		
		\item  If, moreover, there exists a sequence $p_n\in K$ accumulating on $p$ through a local quadrant $Q$ at $p,$ then at least one of the two  branches adjacent to $Q$ is contained in $K$.
	\end{enumerate}
\end{proposition}

{\it Proof: } Without loss of generality, considering an iterate of $f$ if
	necessary, assume $K$ is invariant, $p$ is fixed, and all four branches at $p$ are $f$%
	-invariant. Also, by the Hartman-Grobman theorem, considering a appropriate $%
	(C^0)$ local chart at $p$ we may assume that $f$ is linear in a neighborhood
	of $p:$  the restriction of $f$ to this neighborhood is given by the
	restriction of a linear map to a neighborhood of the origin, where the
	stable/unstable manifolds are the vertical/horizontal axes.
	
	Item 1. First, suppose $K$ avoids $W^u(p)\setminus \{p\}$. Then, if $%
	W_{loc}^u(p)$ denotes a local unstable manifold at $p$, we have that $K$ is
	disjoint from $\Gamma ^u=\overline{W_{loc}^u(p)}\setminus f^{-1}(W_{loc}^u(p))$,
	and therefore there is a neighborhood $U_u$ of the compact set $\Gamma ^u$
	that is disjoint from $K$. From the local linearization one easily verifies
	that $\bigcup_{k\geq 0}f^{-k}(U^u)$ contains a set of the form $B\setminus
	W^s(p)$ where $B$ is a ball around $p$, which we can choose small enough so
	that $K$ intersects the complement of $B$. From the invariance of $K,$ it
	follows that $K\cap B\subset W^s(p)$, and since it must contain $p$ and a
	point outside of $B$, we deduce that $K$ contains a local stable branch at 
	$p$, as desired.
	
	Suppose now that $K$ contains some point $z\in W^u(p)\setminus \{p\}$. We may assume that $z$ lies in the local unstable manifold $W_{loc}^u(p)$,  in particular, in the positive $x$ axis. Assume also that $K$ avoids some point of the unstable
	branch at $p$ containing $z$ (otherwise there is nothing to be done). Then $K
	$ avoids the whole orbit of this point, so we may assume there exists a
	point $x_0$ in $W_{loc}^u(p)\setminus K$ such that the arc $\Gamma ^u$ of $%
	W^u(p)$ joining $x_0$ to $f(x_0)$ contains $z$ and the endpoints of $\Gamma
	^u$ are disjoint from $K$. In local coordinates, one has $\Gamma
	^u=I^u\times \{0\}$ where $I^u$ is an arc in $\mathbb{R}^{+}$. From the fact
	that $K$ is compact, if $\delta >0$ is sufficiently small, the left and
	right edges of $R_\delta =I^u\times [-\delta ,\delta ]$ are disjoint from $K$%
	. This, together with the fact that if we fix $\delta $ small enough $K$
	intersects both $I^u\times \{0\}$ and the complement of $R_\delta $, implies
	that $K$ contains a continuum $K_\delta \subset R_\delta $ intersecting both 
	$I^u\times \{0\}$ and one of the horizontal edges of $R_\delta $. From this
	we also deduce that if $R_\delta ^{-}=I^u\times [-\delta ,0]$ and $R_\delta
	^{+}=I^u\times [0,\delta ]$, either $K\cap R_\delta ^{+}$ contains a
	continuum which intersects the two horizontal edges of $R_\delta ^{+}$, or a
	similar property holds for $R_\delta ^{-}$. In the first case, by a
	straightforward $\lambda $-lemma type argument we deduce that $K$
	accumulates in the local stable branch at $p$ corresponding to the
	positive $y$ axis; and in the second case $K$ accumulates in the opposite
	stable branch. In either case we deduce that $K$ contains a stable
	branch.
	
	Item 2. Suppose $p_n\in K$ converges to $p$ through the first quadrant.
	Inside the linearizing neighborhood of $p\equiv (0,0),$ consider a square of
	the form $Q_1=[0,\delta ]^2$ (clearly contained in the first quadrant) for
	some sufficiently small $\delta >0$ such that $K\cap (Q_1)^c\neq \emptyset .$
	From our hypothesis, we may assume $p_n\in K\cap interior(Q_1).$
	Let $K_n$ be the connected component of $K\cap Q_1$ that contains $p_n.$ As $K$ intersects the complement of $Q_1,$ $K_n$
	intersects $\partial Q_1$ for all integers $n>0.$ As $\{0\}\times [0,\delta ]$
 and $[0,\delta ]\times \{0\}$ are local stable and unstable branches at $p,$ 
if $K_n$ intersects $%
	\{0\}\times [0,\delta ]\cup [0,\delta ]\times \{0\}$ for some $n>0,$ then a
	straightforward $\lambda $-lemma type argument (as in item 1) implies that $K
	$ contains the unstable branch that contains $[0,\delta ]\times \{0\}$ or the
	stable branch that contains $\{0\}\times [0,\delta ].$ As they are both
	adjacent to the first quadrant, the proof is over in this case.
	
	Assume now that all $K_{n^{\prime }s}$ avoid $\{0\}\times [0,\delta ]\cup
	[0,\delta ]\times \{0\}.$ So, for all $n>0,$ $K_n$ intersects $\{\delta
	\}\times [0,\delta ]\cup [0,\delta ]\times \{\delta \}.$ As $p_n$ converges
	to $p,$ there exists a subsequence $n_i\rightarrow \infty $ such that $%
	K_{n_i}\rightarrow \Theta \subset Q_1$ in the Hausdorff topology, as $%
	i\rightarrow \infty .$ It is straightforward to show that $\Theta $ is a
	subcontinuum of $K$ that contains $p$ and intersects $\{\delta \}\times
	[0,\delta ]\cup [0,\delta ]\times \{\delta \}.$ If $\Theta $ does not
	intersect $interior(Q_1),$ then it contains $\{0\}\times [0,\delta ]$ or $[0,\delta ]\times\{0\}$ and we are
	done (from the $f$-invariance of $K).$ If $\Theta $ intersects $%
	interior(Q_1),$ then $f^n(\Theta )$ accumulates on some branch adjacent to
	the first quadrant for $n\rightarrow \infty $ or for $n\rightarrow -\infty $
	(maybe both, for instance when $\Theta \cap \left( \{0\}\times [0,\delta ]\cup [0,\delta
	]\times \{0\}\right) =\{(0,0)\}.$ As $\Theta$ is contained in the $f$-invariant continuum $K$, the proof is over. 
$\Box$
\vskip0.2truecm

From the previous result and the local linearization of periodic saddles, one easily obtains the following
\begin{corollary}\label{coro:contemramo} If $X,Y$ are adjacent branches of a hyperbolic periodic saddle of $f\in \diff^r(M)$ and $K$ is a periodic continuum which accumulates on a point of $X$ locally from the $Y$-side, then $K$ contains one of the branches $X$ or $Y$.
\end{corollary}

A hyperbolic periodic saddle $p=f^n(p)$ always has a continuation, i.e. there are neighborhoods $U$ of $p$ and $\mathcal{U}$ of $f$ in $\diff^r(M)$, and a continuous map $\mathcal{U}\to U$ which maps $g\in \mathcal{U}$ to the unique fixed point $p_g$ of $g^n$ in $U$  (and $p_g$ is a hyperbolic saddle). Moreover, if $X$ is a stable (or unstable) branch of $f$ at $p$, then its corresponding local branch in $U$ also has a continuation, and as a consequence $X_g$ has a continuation, which is the corresponding stable (or unstable) branch of $g$ at $p_g$. Along this section we will use the above notation for continuations of branches and periodic saddles. Let us recall the following result from \cite{pixton}:


\begin{theorem} 
\label{pixresi} If $M$ is a closed surface, there exists a residual subset $\chi^r(M)\subset \diff^r(M)$ such that if $f\in \chi^r(M)$ then all periodic points of $f$ are hyperbolic, the stable and unstable manifolds of any pair of periodic points are transverse (that is, $f$ is Kupka-Smale), and in addition if $X, Y$ are periodic unstable and stable branches, respectively, the maps $$g\mapsto \overline{X_g}, \quad g\mapsto \overline{Y_g},\quad \text{ and } \quad g\mapsto \overline{X_g\cap Y_g}$$ are continuous at $f$.
\end{theorem}

We also recall a basic perturbation result (see lemma 1.6 of \cite{pixton}):

\begin{lemma}[Perturbation lemma] Suppose $p_1,p_2\in M$ are saddle points of $f\in \diff^r(M)$; $Y_1$ is a stable branch at $p_1$ and $X_2$ is an unstable branch at $p_2$. Given $x\in X_2$, a neighborhood $V$ of $x$ and a neighborhood $\mathcal{U}$ of $f$, one has:

\begin{itemize}
	\item[(a)] If $x\in Y_1\cap X_2,$ then there is $g\in \mathcal{U}$ with $g=f$ off $V$ such that $Y_{1g}$ meets $X_{2g}$ transversely at $y$.

	\item[(b)] If $y$ is sufficiently close to $x$ there is $g\in \mathcal{U}$ with $g=f$ off $V$ such that $g(y)=f(x).$ If $f^k(x)\notin V$ for all $k<0,$ then $y\in X_{2g};$ if $y\in Y_1$ and $f^{k}(y)\notin V$ for all $k>0$ then $y\in Y_{1g}.$
	\end{itemize}
\end{lemma}

The next result is essentially theorem 1.7 of \cite{pixton}:
\begin{theorem}\label{Pix1.7} Let $S$ be a surface and let $p$ be a hyperbolic periodic saddle of $f\in \mathrm{Diff}^r(S).$ Also let $X$ and $Y$ be an unstable and a stable branch at $p$, respectively, such that $\overline {\mathcal{O}_f(X)}$ is planar. Suppose that there exists $y$ such that either $y\in \mathcal{O}_f(X)\cap Y$, or $\mathcal{O}_f(X)\cap Y=\emptyset$ and $y\in \overline{\mathcal{O}_f(X)}\cap Y$. Then, given a $C^r$-neighborhood $\mathcal{U}$ of $f$ and a neighborhood $V$ of $y$, there exists $g\in \mathcal{U}$ such the continuations $X_g$ and $Y_g$ of $X$ and $Y$ intersect transversely in $V$.
\end{theorem}

{\it Proof: } This is a direct consequence of theorem 1.7 of \cite{pixton}, which has the same statement with the additional assumption that $S$ is planar. We are only assuming that $\overline{\mathcal{O}_f(X)}$ is planar, but the proof in \cite{pixton} relies only on local perturbations in small neighborhoods of points of $\mathcal{O}_f(X)$, and therefore it applies to our setting (indeed, modifying $S$ outside a neighborhood of $\mathcal{O}_f(X)$, we may assume that $S$ is planar for the purposes of this proof).
$\Box$

\vskip0.2truecm

As a consequence, we have (similar to corollary 2.3 of \cite{pixton}):

\begin{corollary}\label{Pix2.3} Let $X$ and $Y$ be an unstable and a stable branch of a periodic saddle $p$ of $f\in \chi^r(M)$. If $\overline{\mathcal{O}_f(X)}$ is planar and intersects $Y$, then $X$ intersects $Y$ transversely.
\end{corollary}

{\it Proof: } We may modify the surface outside a neighborhood of $\overline{\mathcal{O}_f(X)}$ and assume that $M$ is planar itself. In this case $Y$ changes, but a small local branch $Y_0\subset Y$ does not. From theorem \ref{Pix1.7}, after a small perturbation supported in a neighborhood of $\overline{\mathcal{O}_f(X)}$ one finds a transverse intersection between the continuations of $X$ and $Y$ (which implies a transverse intersection between the continuation of $X$ and $Y_0$). Since this perturbation is localized in a neighborhood of $\overline{\mathcal{O}_f(X)}$, it can be translated to a perturbation in the original surface, with the same consequence. But since $f\in \chi^r(M)$, this implies that there already existed a transverse intersection between $X$ and $Y$.
$\Box$

\vskip0.2truecm

The proof of the following lemma is essentially contained (with minor variations) in the proof of lemma 2.4 of \cite{pixton}. Note that we do not assume that the surface is planar, so the conclusion is weaker.

\begin{lemma}
\label{Pix2.4} Let $X$ and $Y$ be an unstable and a stable branch of a
periodic saddle $p$ for some $f\in \chi ^r(M)$. Suppose there are sequences $%
(y_k)_{k\geq 0}$ of points and $(n_k)_{k\geq 0}$ of positive integers such
that $y_k\to y\in Y$ and $x_k=f^{-n_k}(y_k)\to x\in X$ as $k\to \infty $.
Under these assumptions, $\overline{\mathcal{O}_f(Y)}$ intersects at least
one of the unstable branches at $p$. In addition:

\vskip0.1truecm

a) If $\overline{\mathcal{O}_f(Y)}$ does not intersect $X$, then there
exists a sequence $(\ell _k)_{k\geq 0}$ with $0<\ell _k<n_k$ such that both $%
(f^{-\ell _k}(y_k))_{k\geq 0}$ and $\mathcal{O}_f(Y)$ accumulate on a point $%
x_{*}\in X_{*},$ the other unstable branch at $p,$ locally from the $Y$-side.

\vskip0.1truecm

b) If $x_k\to x$ locally from the $Y$-side, then there exists an unstable
branch at $p,$ a point $x_0$ in that branch (the branch might be $X$ or not)
and a sequence $(\ell _k)_{k\geq 0}$ with $0<\ell _k\leq n_k$ such that both 
$(f^{-\ell _k}(y_k))_{k\geq 0}$ and $\mathcal{O}_f(Y)$ accumulate on $x_0,$
locally from the $Y$-side.
\end{lemma}

{\it Proof: }
Let $N>0$ be the smallest common period of $X,Y$ and of course $p.$ Consider
an arbitrarily small ball $B_\epsilon (y).$ If $\epsilon >0$ is sufficiently
small, then $f^i(B_\epsilon (y))\cap B_\epsilon (y)=\emptyset $ for all $%
1\leq i\leq N$ and $f^j(y)\notin B_\epsilon (y)$ for all $j>0.$ Now
following lemma 9, perturb $f^{-1}$ inside $B_\epsilon (y)$ so that the
resulting diffeomorphism, $g^{-1},$ is $\epsilon $-$C^1$-close to $f^{-1}$
and it satisfies $g^{-1}(y)=f^{-1}(y_k),$ for some arbitrarily large $k.$ We
have two possibilities:

\begin{enumerate}
\item  $f^{-i}(y_k)\notin B_\epsilon (y)$ for all $1\leq i<n_k.$ In this
case, $g^{-n_k}(y)=f^{-n_k}(y_k)$ which is $\epsilon $-close to $x\in X,$ if 
$k$ is large enough.

\item  the above does not happen. In this case, let $i_0$ be the smallest $%
1<i<n_k$ such that $g^{-i}(y)=f^{-i}(y_k)\in B_\epsilon (y).$ 
There exists a constant $T>0$ which depends only on $f$ (obtained from the 
linearization of $f$ near $p$) such that for some
integer $s>0,$ 
$\mathcal{O}_g(Y_g)\ni g^{-i_0+sN}(y)=f^{-i_0+sN}(y_k)$ is 
$T.\epsilon $-close to a fixed fudamental domain $I^u$ of the unstable manifold of $p,$ from the $Y$-side.
\end{enumerate}

If possibility 1 above happens for all arbitrarily small $\epsilon >0$ and
arbitrarilly large values of $k,$ then as $f\in \chi ^r(M),$ $\mathcal{O}%
_f(Y)$ accumulates on $x\in X.$ And moreover, if $x_k\to x$ locally from the 
$Y$-side (the hypothesis in item ${\it b}$ of the lemma), then as $f\in \chi
^r(M),$ $\mathcal{O}_f(Y)$ also accumulates on $x$ locally from the $Y$-side
(because this happens for $g$). 

And if possibility 2 above happens for all $\epsilon >0$ and all arbitrarily
large values of $k,$ then $\mathcal{O}_f(Y)$ and $f^{-i_0(k)+s(k)N}(y_k)$
both accumulate on some point $\overline{x}\in I^u$ from the $Y$-side as $%
k\to \infty $ ($i_0(k)$ and $s(k)$ are sequences which satisfy $1\leq \ell
_k=i_0(k)-s(k)N<n_k$ and $\ell _k\to \infty $ as $k\to \infty $). Clearly $
\overline{x}$ could belong to $X$. This easily concludes the proof of the
lemma. We just have to note that if the hypothesis in item ${\it a}$  holds, then for all $\epsilon >0$ and all arbitrarily large
integers $k,$ we always fall into possibility 2 above. A final remark is that if the hypothesis in item ${\it b}$ holds, then it is possible that either possibilities 1 or 2 occur. $\Box $

\vskip0.2truecm
\begin{corollary}\label{Coro1} 
	Let $X$ and $Y$ be an unstable and a stable branch of a periodic saddle $p$ of $f\in \chi^r(M)$, such that $\overline{\mathcal{O}_f(Y)}$ intersects $X$. If one of the sets $\overline {\mathcal{O}_f(Y)}$ or $\overline{\mathcal{O}_f(X)}$ is planar, then $Y$ intersects $X$ transversely.
\end{corollary}
{\it Proof: } Corollary \ref{Pix2.3} (applied to $f^{-1}$) deals with the case where $\overline{\mathcal{O}_f(Y)}$ is planar. Now assume that $\overline{\mathcal{O}_f(X)}$ is planar.

There is a sequence $(x_k)_{x\geq 0}$ of points of $\mathcal{O}_f(Y)$ converging to a point $x\in X$. From the fact that $Y$ is a stable branch of $p$, after taking a subsequence we may assume  that $f^{n_k}(x_k)\to y\in Y$ for some sequence $(n_k)_{k\geq 0}$ of positive integers. From lemma \ref{Pix2.4} applied to $f^{-1}$ we see that  $\overline{\mathcal{O}_f(X)}$ intersects $Y$ or there is a sequence $(\ell_k)_{k\geq 0}$ with $0<\ell_k<n_k$ such that $f^{\ell_k}(x_k)$ converges to $y_*\in Y_*$ ($Y_*$ is the other stable branch at $p$) and $\mathcal{O}_f(X)$ accumulates on $y_*,$ both from the $X$-side. In the first possibility, corollary \ref{Pix2.3} implies that $X$ intersects $Y$ transversely and the proof is over. And in the second, it implies that $X$ intersects $Y_*$ transversely. Since $f^{\ell_k}(x_k)\in \mathcal{O}_f(Y)$ converges to $y_*$ from the $X$-side, it follows from proposition \ref{Rectangles} that $\mathcal{O}_f(Y)$ intersects $X$. Thus there is $i\geq 0$ such that $f^i(Y)$ intersects $X$. As this intersection is transverse and since $X$ and $Y$ are branches of the same periodic saddle, by a repeated application of the $\lambda$-lemma we conclude that $f^{ki}(Y)$ intersects $X$ transversely, for each $k\geq 1$. In particular since $Y$ is a branch of a periodic saddle, there exists $n$ such that $f^n(Y)=Y$, and choosing $k=n$ we conclude that $Y$ intersects $X$ transversely.
$\Box$

\vskip0.2truecm

\begin{theorem}\label{Pix-new} Let $p,q$ be periodic saddles of $f\in \chi^r(M)$, and suppose there is an unstable branch $Z$ of $q$ and a stable branch $Y$ of $p$ such that $\overline{\mathcal{O}_f(Z)}$ intersects $Y$. If one of the sets $\overline{\mathcal{O}_f(Z)}$ or $\overline{\mathcal{O}_f(Y)}$ is planar, then $\mathcal{O}_f(Z)$ intersects $Y$ transversely.
\end{theorem}

{\it Proof: } Suppose by contradiction that $\mathcal{O}_f(Z)\cap Y=\emptyset$, and let $y\in\overline{\mathcal{O}_f(Z)}\cap Y$, so there exists a sequence $(z_k)_{k\geq 0}$ of points of $\mathcal{O}_f(Z)$ such that $z_k\to y$ as $k\to \infty$.

Suppose first that there exists a neighborhood of $y$ which is avoided by the preorbit of $z_k$ for arbitrarily large values of $k$. Then, the Perturbation lemma allows us to find $g$ arbitrarily close to $f$ such that the continuations of $\mathcal{O}_f(Z)$ and $Y$ for $g$ intersect transversely. The fact that $f\in \chi^r(M)$ then implies that $\mathcal{O}_f(Z)$ already intersected $Y$ transversely prior to the perturbation, contradicting our assumption.

In the remaining case, we can assume that there exists a sequence $(m_k)_{k>0}$ of positive integers such that $f^{-m_k}(z_k)\to y$. Since $z_k\to y$, and $y$ lies in a stable branch of $p$, using the linearization of $f$ at $p$ one easily verifies that, after replacing $(z_k)$ by a subsequence, there are numbers $0<n_k<m_k$ and an unstable branch $X$ of $p$ such that $f^{-n_k}(z_k)$ converges to some point $x\in X$ from the $Y$-side. Lemma \ref{Pix2.4} allows us to conclude that there is an unstable branch $X_*$ of $p$ ($X_*$ could be equal to $X$) such that both $\mathcal{O}_f(Z)$ and $\mathcal{O}_f(Y)$ accumulate on a point $x_*\in X_*$ locally from the $Y$-side. Thus there exist $i,j\in \Z$ such that both $f^i(Y)$ and $f^j(Z)$ accumulate on $x_*$ locally from the $Y$-side. 

Let us prove that $Y$ meets $X_*$ transversely. Recall that one of the sets $\overline{\mathcal{O}_f(Y)}$ or $\overline{\mathcal{O}_f(Z)}$ is planar. 
Suppose first that $\overline{\mathcal{O}_f(Y)}$ is planar. Then since it accumulates on a point of the unstable branch $X_*$ of the same saddle, corollary \ref{Coro1} implies that $X_*$ intersects $Y$ transversely. 
Now suppose that $\overline{\mathcal{O}_f(Z)}$ is planar. As $f^j(Z)$ accumulates on $x_*$ locally from the $Y$-side, corollary \ref{coro:contemramo} implies that $\overline{f^j(Z)}$ contains either $Y$ or $X_*$, and therefore $\overline{\mathcal{O}_f(Z)}$ contains either $\overline{\mathcal{O}_f(Y)}$ or $\overline{\mathcal{O}_f(X_*)}$. Thus one of the last two sets is planar, hence again corollary \ref{Coro1} implies that $X_*$ meets $Y$ transversely.

Since $X_*$ meets $Y$ transversely and $\mathcal{O}_f(Z)$ accumulates on $x_*\in X_*$ locally from the $Y$-side, proposition \ref{Rectangles} implies that $\mathcal{O}_f(Z)$ intersects $Y$ transversely. This contradicts the assumption at the beginning of the proof. $\Box$
\vskip0.2truecm

Recall from the introduction that when $M$ is a closed surface of positive genus, 
$\FM_0^r(M)$ denotes the elements of $\diff^r_0(M)$ which have a contractible 
periodic saddle with a full mesh. If the genus of $M$ is $0$ we simply let 
$\FM_0^r(M)=\diff^r_0(M)$.

\begin{theorem}
	\label{comeconovo} Suppose $M$ is a closed orientable surface  
and $f\in \chi_{FM}=\chi^r(M)\cap \FM^r_0(M)$. If $Y$ is a stable branch of a hyperbolic periodic saddle $p$ and $X$ is an unstable branch of a hyperbolic periodic saddle $q$, such that $\overline{\mathcal{O}_f(X)}$ intersects $Y,$ then there exists $i\in \Z$ such that $Y$ intersects $f^i(X)$ transversely. Moreover, if $p$ and $q$ are in the same orbit, then one may choose $i=0$.
\end{theorem}
{\it Proof: } If $M$ is a sphere, theorem B of \cite{pixton} concludes the proof.
So assume that $M$ has positive genus. 
If one of the sets $\overline{\mathcal{O}%
_f(Y)}$ or $\overline{\mathcal{O}_f(X)}$ is contractible, the proof  follows
form theorem \ref{Pix-new}. So suppose neither of those sets is
contractible. Following proposition \ref{bounddiam} we can choose some
periodic saddle $z^{\prime }$ whose rotation vector is different from the
rotation vectors of $p$ and $q,$ and compact arcs $\lambda _s\subset
W^s(z^{\prime })$ and $\lambda _u\subset W^u(z^{\prime })$ such that $%
z^{\prime }\in \lambda _s\cap \lambda _u$ and $M\setminus (\lambda _s\cup
\lambda _u)$ is a union of open topological disks, whose  covering diameters
are bounded from above by some constant $M(f)>0.$ Let $N>0$ be a common period of 
$X,Y,z^{\prime }$ and all four branches at $z^{\prime }.$ And let $
\widetilde{Y}$ be a connected component of $\pi ^{-1}(Y).$ Assume it is
bounded. Then $closure(\widetilde{Y})$ is a bounded continuum which
satisfies $\widetilde{f}^N(closure(\widetilde{Y}))=g(closure(\widetilde{Y})),
$ where $\widetilde{f}$ is the lift of $f$ to the universal covering used to
compute the rotation set, and $g$ is a Deck transformation. As the rotation
vector of $p$ (which is equal to $[g]/N$) is different from the rotation
vector of $z^{\prime },$ $closure(\widetilde{Y})$ avoids $\pi ^{-1}(\lambda
_s\cup \lambda _u).$ So $\overline{Y}$ avoids $\lambda _s\cup \lambda _u.$
The same happens for $f^i(\overline{Y})$ for all $0\leq i\leq N-1.$ As $
\overline{\mathcal{O}_f(Y)}=\overline{Y}\cup f(\overline{Y})\cup ...\cup
f^{N-1}(\overline{Y}),$ we get that $\overline{\mathcal{O}_f(Y)}$ is
contractible, a contradiction with our initial assumption. So $\widetilde{Y}$
is unbounded, therefore it intersects $\pi ^{-1}(\lambda _s\cup \lambda _u)$
and finally, we get that $Y$ must have a transversal intersection with $%
\lambda _u.$

Analogously, $X$ must have a transversal intersection with $\lambda _s.$ As $%
\lambda _s$ and $\lambda _u$ are adjacent branches of the same saddle, the $\lambda 
$-lemma implies that $Y$ intersects $X$ transversely, concluding the proof.
$\Box$
\vskip0.2truecm

{\bf Remark: } Under the notation of the previous theorem, a simple and useful 
consequence is the following: if $q$ has a full mesh and $p$ is a 
hyperbolic periodic saddle point contained in $\overline{W^u(q)},$ then 
$W^u(q)$ has a
transversal intersection with $W^s(p).$ If $p=q$, this is trivial because
$q$ has homoclinic intersections (it has a full mesh!). And if $q \neq p$, then as 
$p \in \overline{W^u(q)}$, some unstable branch at $q$ accumulates on some point 
of a stable branch at $p,$ so the theorem can be applied in order to obtain the 
intersection. We just have to note that 
$\overline{W^u(q)}=\overline{W^u(f^i(q))},$ for all integers $i$ (see observation 1
below).

\vskip0.2truecm

\section{Main technical result}

In this section we prove a result which contains the proofs of theorem 1 and
theorem D.
Recall that the set $\chi_{FM}\subset \diff_0 ^r(M)$
was defined in the statement of theorem \ref{comeconovo} and for any 
given subset $K \subset M$, 
$Filled(K)$ is the union of $K$ with all the inessential connected components
of the complement of $K.$ 

Before stating the theorem, we make two observations:

\begin{description}
\item[Observations.]   {\it Let $M$ be a closed orientable surface of
genus $g$ and $z\in M$ be a hyperbolic $n$-periodic saddle point for 
some $f\in \chi_{FM}\subset \diff_0 ^r(M).$

\begin{enumerate}
\item  If $g>0$ and $z$ has a full mesh, then for all $0\leq i\leq n-1$ as $W^{s\text{ }%
(or\text{ }u)}(z)$ has transversal intersections with $W^{u\text{ }(or\text{
}s)}(f^i(z)),$ we get that $\overline{W^s(z)}=\overline{W^s(f^i(z))}$ and $
\overline{W^u(z)}=\overline{W^u(f^i(z))};$

\item  For any $g \geq 0$, if $\overline{W^s(z)}$ is contractible, 
then the complement of 
$Filled(\cup _{i=0}^{n-1}f^i\overline{(W^s(z)}))$ is an $f$-invariant fully
essential open set $D\subset M$ (this follows from proposition \ref{bounddiam}).
So, $\partial D$ has finitely many (contractible) connected components 
$K_0,K_1,...,K_{l-1}\subset \cup _{i=0}^{n-1}f^i(\overline{W^s(z)})$ 
for some $1\leq l\leq n,$
such that $f(K_i)=K_{i+1}$ (for $i=0,1,...,l-2)$ and $f(K_{l-1})=K_0.$ As 
$f^n(K_i)=K_i$, the prime ends rotation number of $f^n\mid_{K_i}$ is the same 
for all $0\leq i \leq l-1.$ So it makes sense to say that the prime ends 
compactification of $D$ has rational (or irrational) rotation number: the 
rotation number is well-defined at each boundary component of $D$ and it assumes
the same value at all of them 
(see \cite{koropatmey} for precise
definitions of prime ends in this kind of situation). Clearly, analogous properties 
hold for any periodic inessential connected component of the complement of 
$\cup _{i=0}^{n-1}f^i\overline{(W^s(z)}).$ 
\end{enumerate}}

\end{description}

\begin{theorem}
\label{repulsor0} Let $M$ be a closed orientable surface,  
$f\in \chi _{FM}\subset \diff_0^r(M),$ for any $r\geq 1,$ 
and $z\in M$ be a hyperbolic $n$-periodic saddle point which satisfies:

\begin{enumerate}
\item  Either, the genus of $M$ is positive and $z$ has a full mesh, or 
$\overline{W^s(z)}$ is contractible (for any $genus \geq 0$);
%

\item  When $z$ has a full mesh, there exists a $f$-periodic connected
component $D$ of the complement of $\overline{W^s(z)}\subset M$ which has
rational prime ends rotation number and in the second case, as explained in
 observation 2 above, $D$ is a $k$-periodic connected component of 
$(\cup _{i=0}^{n-1}f^i\overline{(W^s(z)})^c$ (maybe fully essential) for which 
the prime ends rotation number of $f^{n.k}$ at each boundary component of $D$
is the same rational number;
\end{enumerate}

Under the above hypotheses, assuming that some integer $N\geq 1$ is a common
period of $z,$ $D$ and each connected component of $\partial D,$ as well
as all 4 branches at $z$, and also assuming that the prime ends rotation
number of each boundary component of $D$ is equal to $p/q$ (a rational
number), the following conclusions hold:

A) If $V$ is a cross-section in $D$ which satisfies $closure_D(f^{N.q}(V))%
\subset V,$ then $f^{N.q}$ has a fixed point in $V.$

B) Each connected component of the boundary of $\widehat{D},$ the prime ends
compactification of $D,$ has finitely many periodic prime ends, which can
only be of two types: sources and accessible hyperbolic saddles. Moreover,
for each saddle in the later case, its stable manifold avoids $D,$ 
one unstable 
branch is contained in $D$ and thus the whole stable manifold is accessible from $D$.

%
\end{theorem}

{\it Proof: }In both cases above (either when $z$ has a full mesh or when $
\overline{W^s(z)}$ is contractible) we denote by $g=f^N$ the homeomorphism
acting on each connected component of $\partial $$D.$

First, let us show that the map $\widehat{g}:\widehat{D}\rightarrow \widehat{%
D}$ induced by $g,$ does not have intervals of $q$-periodic points in each circle
in $\partial \widehat{D}.$ Indeed, if it had, then there would be a crosscut 
$\alpha $ in $D$ such that all points in $\partial D$ that belong to the
boundary of a cross-section associated to $\alpha $ would be $g^q$-fixed
(because accessible points had to be $g^q$-fixed and they are dense in $\partial D)
$, a contradiction with the fact that $f$ is Kupka-Smale.

{\bf Proof of A.}{\it \ }Assume $C$ is a crosscut in $D$ such that $%
g^q(C)\cap C=\emptyset $ and $V$ is a cross-section associated to $C$ that
satisfies $g^q(V)\subset V.$ Consider an arc $\gamma $ contained in the
region between $C$ and $g^q(C),$ whose endpoints are, some $a\in C$ and $%
g^q(a)\in g^q(C).$ Clearly, $g^{q.n}(\gamma )\cap \gamma =\emptyset $ for
all $n\notin \{-1,0,1\}.$ Let us look at $\cup _{m\geq 0}g$$^{q.m}(\gamma )$
and the set of accumulation points of $g^{q.m}(\gamma )$ as $m\rightarrow
\infty ,$ which we call the $\omega $-limit set of 
$\cup _{m\geq 0}g^{q.m}(\gamma ).$ It clearly is a $g^q$-invariant continuum 
$K\subset (V\cup \partial
V)\backslash C.$ From the hypothesis on $C,$ we get that $K\cap \left( \cup
_{m\geq 0}g^{q.m}(\gamma )\right) =\emptyset .$

Even when $D$ is homotopically 
unbounded, proposition \ref{geral} together with proposition \ref{bounddiam}
imply that $K$ is contractible. Thus the main theorem of \cite{andresrata}
states that $K$ is a rotational attractor, or a rotational repeller, or it 
contains
a $g^q$-fixed point. When $K$ is a rotational attractor (or repeller), we mean 
that it is a contractible continuum
in $M$ which may, or may not, separate $M,$ and it attracts (or repells) all 
nearby points that belong to at least one connected component of its complement. 
In case $K$ separates $M,$ one of the connected components of $K^c$ is fully 
essential and the others are open disks. Clearly, by the Cartright-Littlewood theorem 
\cite{cartlit}, $Filled(K)$ always contains $g^q$-fixed points, but when $K$
separates $M$, the $g^q$-fixed points may not belong to $K.$ Note that this
definition is slightly different from the one which appears in \cite
{andresrata}.


So, when $K\subset V,$ as $V$ is a topological open disk, if $K$ is a
rotational attractor or repeller, then there is a $g^q$-fixed point in $%
Filled(K)\subset V.$ Thus, from the three possibilities given by the main
theorem of \cite{andresrata}, when $K\subset V,$ $V$ always contains a $g^q$%
-fixed point (as a matter of fact, in our case $K$ can not be a repeller
because it attracts positive iterates of $\gamma $).

Thus, let us assume that $K$ intersects $\partial D.$

First, we deal with the case when $K$ is a rotational attractor. It implies that 
$z\in K.$ When $z$ has a full mesh, as $W^u(z)$ is unbounded, this is not 
possible. And when $\overline{W^s(z)}$ is
contractible, it would imply that the whole $W^u(z)$ is contained in $K$, 
therefore $\overline{W^u(z)}$ is also contractible. Let us keep this information
for a while.

As $K$ can not be a rotational repeller, the remaining possibility
given by the main result of \cite{andresrata} is that $K$ must have $g^q$-fixed
points. If at least one of these points belongs to $V,$ then we are
done.

So, assume all of them belong to $\partial D.$ As $f$ only has hyperbolic
periodic points and positive iterates of $\gamma $ accumulate on $K,$ the
$g^q$-fixed points in $K$ can not be sources. They can not be sinks either,
because 

\begin{equation}
\label{melhorapoqq}\partial D\subseteq \cup _{i=0}^{n-1}f^i\overline{(W^s(z)}%
).
\end{equation}

So the $g^q$-fixed points in $K$ are all saddles. More precisely:

\begin{proposition}
\label{whatsapp} $K$ contains a hyperbolic $g^q$-fixed saddle point of index 
$-1.$
\end{proposition}

{\it Proof:} If $K$ is a rotational attractor, $z$ belongs to $K$ and 
we are done. So, we can
assume that $K$ is not a rotational attractor.

By contradiction, suppose that all $g^q$-fixed points in $K$ are
orientation-reversing saddles, denoted $\{w_1,w_2,...,w_m\},$ that is, they
all have index $1.$ As there are finitely many such points in $K$ (because $f$
 is Kupka-Smale), we can blow up each $z_i$ into a circle $C_i,$ such that
the induced dynamics (by $g^q$) on $C_i$ is that of a semi-rotation. In the
disk $D_i$ bounded by $C_i,$ we define a dynamics in the natural way: the
extended map fixes the center of each $D_i$ and semi-rotates each concentric
circle. In this new space obtained after blowing up all $g^q$-fixed points in 
$K$ and adding disks, positive iterates of $\gamma $ still accumulate on an 
invariant continuum, which can not be a rotational attractor, because
prior to the blow-ups we assumed it was not. The problem is that this
continuum has no $g^q$-fixed points. This contradicts the main result of 
\cite{andresrata} and proves the proposition. $\Box $

\vskip0.2truecm

So $K\cap \partial D$ contains at least one hyperbolic  
$g^q$-fixed saddle point $w$ of index $-1$ (note that $w$ can be equal to $z$,
for instance when $K$ is a rotational attractor).
 From expression (\ref{melhorapoqq}), if $w$ is not equal to $z,$ 
for some $0\leq i\leq n-1,$ $f^i(W^s(z))=W^s(f^i(z))$ accumulates on $w$
and therefore, by theorem \ref{comeconovo}, 
$W^s(f^j(z))$ has a $C^1$%
-transverse intersection with $W^u(w)$ for some $0\leq j\leq n-1$ $.$ As $D$
is a connected component of the complement of $\cup _{i=0}^{n-1}f^i\overline{%
(W^s(z)}),$ we finally get that $\overline{W^s(w)}\cap D=\emptyset .$ Clearly,
the last intersection is also empty when $w=z$, by the choice of $D.$

Whether $K$ is a rotational attractor or not, 
as $K$ is a $g^q$-invariant continuum that contains $w,$ it contains at
least one of the four branches at $w,$ see item 1 of proposition \ref
{contemramo}. Let $\lambda $ be a branch at $w,$ either stable or unstable,
contained in $K.$ We know that $\cup _{m\geq 0}$$g^{q.m}(\gamma )$
accumulates on the whole $\lambda .$ In the following we will show that this
accumulation implies that $\cup _{m\geq 0}g$$^{q.m}(\gamma )$ has a
topologically transverse intersection with $W^s(w^{\prime }),$ for some
saddle $w^{\prime }\in \partial D.$ As we explained above, for any such a 
saddle $w^{\prime},$  
$\overline{W^s(w^{\prime })}\cap D=\emptyset ,$ a contradiction because $\cup
_{m\geq 0}g$$^{q.m}(\gamma )$ is contained in $D.$

The fact that $\lambda $$\subset K,$ implies that $\overline{\lambda }$ is
contractible, so again from the main theorem of \cite{andresrata} we get
that $\Theta ,$ the accumulation set of $\lambda $ ($\Theta $ is given by
the $\omega $ or $\alpha $-limit set of $\lambda ,$ depending on whether $%
\lambda $ is an unstable or stable branch), is a rotational attractor
or it contains a $g^q$-fixed point (it can not be a repeller because it is
accumulated by positive iterates of $\gamma ).$ If $\Theta $ is a rotational
attractor, then it either belongs to $V$ and we are done (because in this
case, $Filled(\Theta )\subset V$ and it contains $g^q$-fixed points), or it
intersects $\overline{W^s(z)}.$ In this last case, $z$ must belong to $%
\Theta.$ So, either part A of the theorem is proved or all $g^q$-fixed points 
in $\Theta $ belong to 
$\partial D.$ In particular, some point $w_1\in \Theta \cap \partial D$ 
($w_1$ may be equal to $z$) is a 
$g^q$-fixed hyperbolic saddle of index $-1$ (again, this follows from
proposition \ref{whatsapp} and the fact that a $g^q$-fixed point in
$\Theta \cap \partial D$ can not be a source or a
sink since it is accumulated by positive iterates of $\gamma $ and by $W^s(z)
$). As $\Theta \subset K$ is a $g^q$-invariant continuum, item 1 of
proposition \ref{contemramo} implies that $\Theta $ contains some branch at $%
w_1.$ This branch has a contractible closure and as before, it accumulates
on a continuum $\Theta _2.$ Exactly as we did above, if $\Theta _2$ is a 
rotational attractor, then 
$z$ belongs $\Theta _2.$ So, in any case (see proposition \ref{whatsapp}),
 there is a $g^q$-fixed 
hyperbolic saddle of index $-1,$ denoted 
$w_2\in \Theta _2\cap \partial D.$ Clearly, $\Theta _2$ contains a branch at 
$w_2$ and so, $\lambda \stackrel{def}{=}\lambda _0$ accumulates on a whole
branch $\lambda _1$ at $w_1$ and this branch at $w_1$ accumulates on a whole
branch $\lambda _2$ at $w_2$ and so on. Now the proof goes as follows: 
As $K$ is compact and $g^q$ 
is Kupka-Smale, some
point $w_N$ (for some $N\geq 0$) in the sequence $w\stackrel{def}{=}%
w_0\rightarrow w_1\rightarrow w_2\rightarrow ...\rightarrow w_N\rightarrow
w_{N+1}\rightarrow ...\rightarrow w_N$ must appear twice. Note that each $w_i
$ is a $g^q$-fixed hyperbolic saddle of index $-1$ contained in $K\cap
\partial D$ and $w_i\rightarrow w_{i+1}$ means that the branch $\lambda
_i\subset K$ at $w_i$ accumulates on the whole branch $\lambda _{i+1}$ at $%
w_{i+1}.$

As some contractible branch $\lambda _N$ at $w_N$ accumulates on $w_N,$
theorem \ref{comeconovo} implies that $\lambda _N$ has a transversal
homoclinic point. As $\lambda _N\subset K$ and $\cup _{m\geq 0}$$%
g^{q.m}(\gamma )$ accumulates on the whole $\lambda _N,$ proposition
\ref{Rectangles} implies that $\cup _{m\geq 0}g$$^{q.m}(\gamma )$ must have
a topologically transverse intersection with $W^s(w_N)$ (because it has to
enter the horseshoe rectangle). To see this, assume, by contradiction, that
positive iterates of $\gamma $ under $g^q$ only intersect $W^u(w_N).$ This
would imply that for infinitely many large integers $m>0,$ $g^{q.m}(\gamma )$
intersects the boundary of a fixed horseshoe rectangle at $w_N,$ only at
points belonging to $W^u(w_N).$ So, $w_N\in \gamma ,$ a contradiction
because $\gamma \subset D$ and $w_N\in \partial D.$

At last, as $W^s(w_N)\subset \overline{W^s(f^j(z)),}$ for some $0\leq j\leq
N-1,$ we get that $\overline{W^s(w_N)}\cap D=\emptyset .$ As $\cup _{m\geq 0}
$$g^{q.m}(\gamma )$ is contained in $D,$ this is the final contradiction
that proves A (see figure 1).

\begin{figure}[!h]
	\centering
	\includegraphics[scale=0.43]{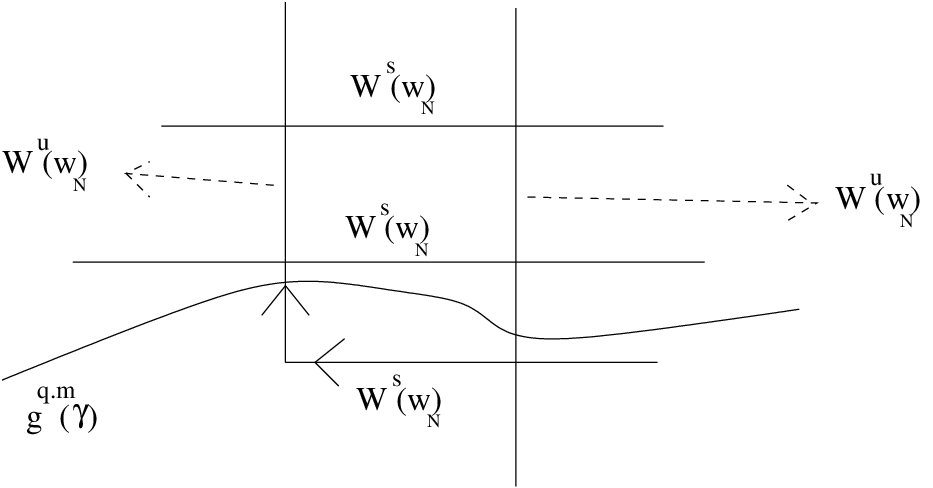}
	\caption{Diagram showing iterates of $\gamma$ entering a rectangle of a
          horseshoe.}
	\label{clopen1}
\end{figure}

\vskip 0.2truecm


{\bf Proof of B.}{\it \ }Now we show that there are only finitely many
periodic prime ends in each of the (finite) connected components of $%
\partial \widehat{D}.$ Suppose not. Then, the complement of the $\widehat{g}^q$
-fixed point set in some component of $\partial \widehat{D}$ has infinitely many
connected components, which are disjoint open intervals. In each of these
intervals, orbits under $\widehat{g}^q$ move in one direction, so for each 
interval,
one of its boundary points is $\widehat{g}^q$-fixed and attracts orbits from at least one
side.

The important conclusion is that the existence of infinitely many 
$\widehat{g}^q$-fixed prime ends imply the existence of infinitely many 
$\widehat{g}^q$-fixed prime
ends that attract from at least one side. Let $\widehat{p}^{*}\in \partial 
\widehat{D}$ be a $\widehat{g}^q$-fixed prime end that attracts from at least one
side. If $(C_n^{*})_{n\geq 0}$ is a prime chain that represents $\widehat{p}%
^{*},$ by choosing diameter of $C_0^{*}$ sufficiently small, we can assume
that there are no $g^q$-fixed points in the cross-section $V_0^{*}.$ This
follows from the fact that $f:M\rightarrow M$ has finitely many $n$-periodic
points for each integer $n\geq 1.$

As $\widehat{p}^{*}$ attracts from at least one side, part A above implies
that for all $n\geq 0,$ {\bf $g^q(C_n^{*})\cap C_n^{*}\neq \emptyset .$} So,
all principal points associated to $\widehat{p}^{*}$ are $g^q$-fixed. As the
principal set is connected, there is only one such point, denoted $p^{*}\in
\partial D$ and $\widehat{p}^{*}$ is an accessible prime end. The $g^q$%
-fixed point $p^{*}$ can not be a sink because it belongs to the closure of
stable manifolds. Neither a source, because $\widehat{p}^{*}$ is accessible
and attracts from at least one side, see lemma 6.3 of \cite{yorke}. So it is
a saddle. As $p^{*}\in \partial D\subseteq \cup _{i=0}^{n-1}f^i\overline{%
(W^s(z)}),$ as we did in part A above, if $%
p^{*}\notin Orb(z),$ then at least one branch of $W^u(p^{*})$ has a
transverse intersection with $W^s(f^i(z))$ for some $0\leq i\leq n-1.$ So, $
\overline{W^s(p^{*})}\subset \overline{W^s(f^i(z))}$ which implies that $
\overline{W^s(p^{*})}\cap D=\emptyset .$ Nevertheless, even when $p^{*}\in
Orb(z),$ $\overline{W^s(p^{*})}\subset \cup _{i=0}^{n-1}f^i(\overline{W^s(z)}%
).$

As $p^{*}$ is accessible, there exists an end-cut $\eta :[0,1]\rightarrow 
\overline{D}$ such that $\eta $$([0,1[)\subset D$ and $\eta $$(1)=p^{*}.$
Clearly, considering small neighborhoods of $p^{*},$ we can assume that $%
\eta $ is as in figure 2, essentially because $D$ avoids $\overline{%
W^s(p^{*})}.$

\begin{figure}[!h]
	\centering
	\includegraphics[scale=0.43]{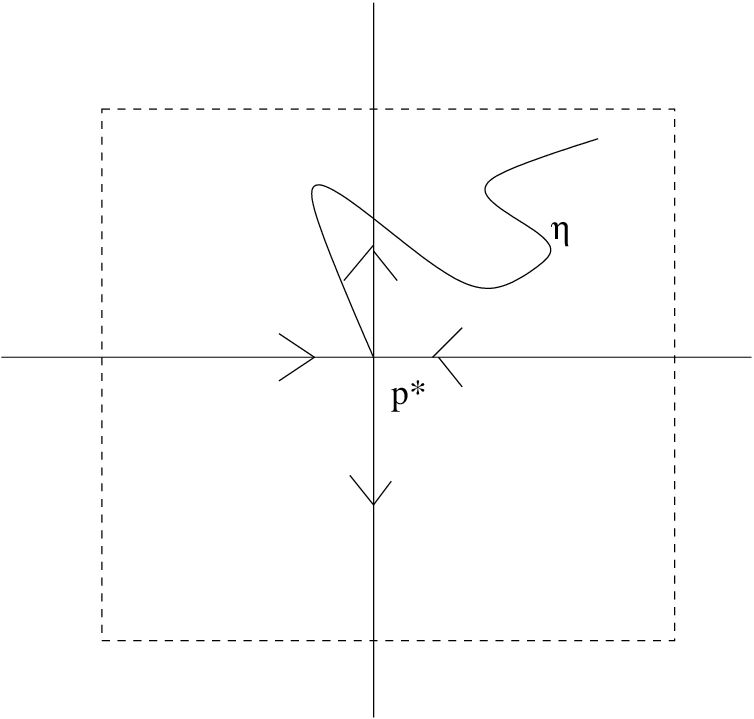}
	\caption{Diagram showing the end-cut $\eta.$}
	\label{clopen2}
\end{figure}

Let $\lambda _{p^{*}}^u$ be the branch of $W^u(p^{*})$ which does not have
an (transverse, as $f\in \chi^r(M) )$ intersection with $\cup
_{i=0}^{n-1}f^i(W^s(z)).$ There must exist such a branch, otherwise $p^{*}$
would not be accessible. In case $p_{*}\notin Orb(z),$ we know that one
unstable branch at $p^{*}$ intersects $\cup _{i=0}^{n-1}f^i(W^s(z))$ and the
other not. So, each of them is $g^q$-invariant and therefore the index of $%
p_{*}$ under $g^q$ is $-1.$ Clearly, as $z$ and all four branches at $z$ are 
$g$-invariant, the index of $z$ under $g^q$ is also $-1.$

From proposition \ref{arcoevita}, the endcut $\eta $ can be deformed into
one that does not intersect $\lambda _{p^{*}}^u$ (in other words, $\eta $ is
homotopic in $D$ to an endcut which avoids all branches at $p^{*}).$ So, we
can assume that $\eta $ is contained in one of the local quadrants adjacent
to $\lambda _{p^{*}}^u$ and thus, some branch of $W^s(p^{*})$ is made of
accessible points in $\partial $$D,$ see proposition \ref{metadeacess}.
These two results appear in \cite{yorke} and below we present proofs of both.


\begin{proposition}
\label{arcoevita} The arc $\eta $ associated to $\widehat{p}^{*}$ can be
chosen so that $\eta \cap \lambda _{p^{*}}^u=\emptyset .$
\end{proposition}

{\it Proof: }By contradiction, assume the proposition does not hold. Then,
there must be a prime chain $(C_n^{*})_{n\geq 0}$ representing $\widehat{p}%
^{*},$ such that each crosscut $C_n^{*}$ is contained in $\lambda _{p^{*}}^u.
$ So, $g^q(C_n^{*})\cap C_n^{*}=\emptyset $ for all $n\geq 0$ and if 
diameter of $C_0^{*}$ is sufficiently small, then by part A of this
theorem, $\widehat{p}^{*}$ does not attract from any side. This
contradiction proves the proposition. $\Box $

\vskip 0.2truecm

\begin{proposition}
\label{metadeacess} There exists a branch of $W^s(p^{*})$ made of accessible
points in $\partial $$D.$
\end{proposition}

{\it Proof: }This proof appears in \cite{frankslecal} in detail and
originally in Mather's paper \cite{mather}. We shortly present it here.

By the Hartman-Grobman theorem, assume we are inside an open ball centered
at $p^{*}$ where the dynamics of $g^q$ is that of $(x,y)\rightarrow (x/2,2y)$
(clearly $p^{*}$ is identified with $(0,0)$). Also assume that $\eta $
belongs to the first quadrant, so $\lambda _{p^{*}}^u$ contains the local
positive $y$-axis. For each $t\in \eta ,$ let $V_t$ be a vertical segment
starting at $t$ and going downwards, until it touches $\partial D.$ As $
\overline{W^s(p^{*})}\cap D=\emptyset ,$ $V_t$ does not cross the $x$-axis.
Analogously, let $H_t$ be an horizontal segment starting at $t$ and going to
the left, until it touches $\partial $$D.$

$\bullet $ If for some $t\in \eta ,$ $H_t$ intersects the $y$-axis, then
there exists a simple closed curve $\beta $ which is made of a vertical arc
in $\lambda _{p^{*}}^u$ from $p^{*}$ to a point in $H_t,$ an horizontal arc
from that point to some point in $\eta $ and a sub-arc of $\eta ,$ going
back to $p^{*}$. Let $B$ be the interior of $\beta .$ If, by contradiction, $%
\cup _{i=0}^{n-1}f^i(W^s(z))$ accumulates on $p^{*}$ through the first
quadrant, then for some $0\leq i^{*}\leq n-1,$ $f^{i^{*}}(W^s(z))$
intersects $B.$ From the choices of $\eta $ and $H_t,$ and the fact that $%
f^{i^{*}}(W^s(z))$ is connected and accumulates on the whole $W^s(f^j(p^{*}))
$ for some integer $j,$ and this stable manifold is not contained in $B,$ $%
f^{i^{*}}(W^s(z))$ must have a transversal intersection with $\lambda
_{p^{*}}^u.$ As we already said, this implies that $p^{*}$ is not
accessible. So $\cup _{i=0}^{n-1}f^i(W^s(z))$ does not accumulate on any
point of the stable branch at $p^{*}$ which locally coincides with the
positive $x$-axis, through the first quadrant. In other words, the whole
local first quadrant is contained in $D.$

\vskip0.2truecm

$\bullet $ Now assume that for all $t\in \eta ,$ $H_t$ is contained in the
first quadrant. Clearly, $V_t\cup H_t$ is a crosscut associated to $\widehat{%
p}^{*}.$ As $t$ gets closer to $p^{*},$ $diameter(V_t\cup H_t)$ becomes
arbitrarily small. So, as $\widehat{p}^{*}$ attracts from at least one side,
for all $t$ sufficiently close to $p^{*},$ part A above implies that, 
$$
g^q(V_t\cup H_t)\cap (V_t\cup H_t)\neq \emptyset . 
$$

But the dynamics near $p^{*}$ is that of $(x,y)\rightarrow (x/2,2y),$ so $%
g^q(V_t)$ is a vertical segment and $g^q(H_t)$ is an horizontal one. This
implies that $g^q(V_t)\cap (H_t)\neq \emptyset .$ In other words, there are
points $s\in int(V_t)$ and $g^q(s)\in int(H_t)$ such that the arc $\zeta
_s\subset V_t\cup H_t$ whose end points are $s$ and $g^q(s),$ satisfies $%
g^q(\zeta _s)\cap \zeta _s=g^q(s).$ if we perform this construction for all $%
t\in \eta $ sufficiently close to $p^{*},$ we obtain a fundamental domain of
the first quadrant entirely contained in $D.$ So the whole
local first quadrant is contained in $D$ and thus the stable branch at $p^{*}$
which locally coincides with the positive $x$-axis is made of accessible
points through the local first quadrant. $\Box $ %

\vskip 0.2truecm

Let us show now that $\lambda _{p^{*}}^u$ is contained in $D$. As all its points 
are accessible from $D$, if it is not contained in $D$, then
$$
\lambda _{p^{*}}^u\cap (\cup_{i=0}^{n-1}f^i(\overline{W^s(z)}))\neq\emptyset.  
$$
 So, theorem \ref{comeconovo} implies that $\lambda _{p^{*}}^u$ would have a 
transversal intersection with 
$f^i(W^s(z)),$ for some integer $i$, a contradiction with the accessibility of $p^*$. 
And finally, as a local quadrant adjacent to $\lambda _{p^{*}}^u$ belongs 
to $D$, the whole $\lambda _{p^{*}}^u$ is contained in $D,$ so 
$W^s(p^*)$ is contained in $\partial D$ and is made of accessible points.

If $\widehat{p}_1^{*},\widehat{p}_2^{*}\in \partial 
\widehat{D}$ are $\widehat{g}^q$-fixed prime ends that attract from at least
one side, then they correspond to accessible saddles in $\partial D.$ Assume
they both correspond to the same saddle $p_{*}\in \partial D.$

If $p^{*}\in Orb(z)$ and $\overline{W^s(z)}\cap W^u(z)=\emptyset $ (this is
only possible when $\overline{W^s(z)}$ is contractible), then there might be
two different $\widehat{g}^q$-fixed prime ends $\widehat{p}_1^{*},\widehat{p}_2^{*}$ which
correspond to $p^{*};$ \ $\widehat{p}_1^{*}$ is associated to one unstable
branch at $p^{*}$ and $\widehat{p}_2^{*}$ is associated to the other.
Clearly, in general, no more than two different prime ends are associated to
a point in $Orb(z).$

Now, assume $p^{*}\notin Orb(z).$ We know that $\overline{W^s(p^{*})}\subset
\cup _{i=0}^{n-1}f^i(\overline{W^s(z)})$ and there are two end-cuts $\eta
_1,\eta _2:[0,1]\rightarrow \overline{D}$ such that $\eta _{1,2}$$%
([0,1[)\subset D$ and $\eta _{1,2}$$(1)=p^{*}.$ As $p^{*}\notin Orb(z),$ by
proposition \ref{arcoevita} we can assume that each of these endcuts is
contained in one of the local quadrants adjacent to $\lambda _{p^{*}}^u,$
the unstable branch at $p^{*}$ which is contained in $D.$ 
This happens because, from theorem \ref{comeconovo},
 the other unstable branch
at $p^{*}$ intersects $\cup _{i=0}^{n-1}f^i(W^s(z)).$

If $\eta _1$ and $\eta _2$ belong to the same quadrant adjacent to $\lambda
_{p^{*}}^u,$ then they are homotopic in $D$ and thus $\widehat{p}_1^{*}= 
\widehat{p}_2^{*}.$ If they belong to different adjacent quadrants, then
as $\lambda _{p^{*}}^u$ is an endcut at $%
p_{*}$ and $\eta _1$ and $\eta _2$ are both homotopic in $D$ to $\lambda
_{p^{*}}^u,$ they are homotopic to each other. So, in this case we also have $\widehat{p}%
_1^{*}=\widehat{p}_2^{*}.$ As $f:M\rightarrow M$ has finitely many periodic
points for each period, in the first case (when $z$ has a full mesh), $
\overline{W^s(z)}$ has finitely many $g^q$-fixed points and in the second
case $\cup _{i=0}^{n-1}f^i(\overline{W^s(z)})$ also has
finitely many $g^q$-fixed points. So, there must be finitely many periodic
prime ends that attract from at least one side, thus finitely many periodic
prime ends in $\partial \widehat{D}.$

Summarizing, if such a periodic prime end attracts from at least one side,
then it attracts from both sides, it is accessible and it corresponds to a saddle in $\partial D,$ for
which one unstable branch is contained in $D$.

Suppose now that $\widehat{p}\in \partial \widehat{D}$ is an isolated 
$\widehat{g}^q$%
-fixed prime end which repels from both sides and let $(C_n)_{n\geq 0}$ be
any prime chain that represents $\widehat{p}.$ Denote by $(V_n)_{n\geq 0}$
the cross-sections associated to the $C_{n^{\prime }s}.$ As before, by
choosing $C_0$ with a sufficiently small diameter, we can assume that there
is no $g^q$-fixed point in $V_0.$

Either one of the following possibilities hold:

\begin{enumerate}
\item  For all $(C_n)_{n\geq 0}$ as above, there is an infinite sequence 
{\bf $C_{n_i}$ }such that {\bf $C_{n_i}\cap g^q(C_{n_i})\neq \emptyset ;$}

\item  the above does not hold;
\end{enumerate}

In the first possibility, all principal points associated to $\widehat{p}$
are $g^q$-fixed, therefore, from our genericity condition, there is only one
such point $p\in \partial D,$ which is an accessible $g^q$-fixed point.

It can not be a sink. If it were a source, then by lemma 6.3 of \cite{yorke}%
, we could choose some prime chain $(C_n)_{n\geq 0}$ representing $\widehat{p%
}$ for which {\bf $C_n\cap g^q(C_n)=\emptyset $} for all $n\geq 0,$ a
contradiction with our hypotheses. So $p$ is a saddle. And as we explained
above, the stable manifold of $p$ is made of accessible points in $\partial D.$
Therefore $\widehat{p}$ attracts from both sides, a contradiction
with the assumption on $\widehat{p}.$

So the first possibility can not happen, and in the second, there is a prime
chain $(C_n)_{n\geq 0}$ representing $\widehat{p}$ such that {\bf $C_n\cap
g^q(C_n)=\emptyset $} for all $n\geq 0.$ 
Again, there are two sub-cases to consider:

\begin{itemize}
\item  For some $i_0,$ $\widehat{g}^{-q.m}(\widehat{C_{i_0}})\rightarrow \widehat{p}$
as $m\rightarrow \infty .$

\item  The above does not happen.
\end{itemize}

In the first case above, $\widehat{p}\in \widehat{D}$ is a source. And in
the second case, 
there exist $C_{i_0}$ and $C_{i_1},$ $i_0<i_1$ arbitrarily large, $%
C_{i_1}\subset V_{i_0}$ such that $g^{-q}(C_{i_0})\subset V_{i_0},$ $%
g^{-q}(C_{i_0})\cap V_{i_1}=\emptyset $ and $g^{-q.n}(V_{i_0})\cap
C_{i_1}\neq \emptyset $ for all $n\geq 0.$ 
So, $\cap _{n\geq 0}\widehat{g}^{-q.n}(\widehat{V}_{i_0})$ is a closed 
connected $\widehat{g}^q$-invariant subset 
of $\widehat{D}$ that intersects $\widehat{C}%
_{i_1}$ and 
$$
\left( \cap _{n\geq 0}\widehat{g}^{-q.n}(\widehat{V}_{i_0})\right) \cap
\partial \widehat{D}=\widehat{p}. 
$$

As $g^{-q}(C_{i_0})\cap V_{i_1}=\emptyset $ , it follows that $%
g^q(C_{i_1})\subset V_{i_0}.$ So, 
\begin{equation}
\label{invset}\cap _{i\geq 0}g^{q.i}\left[ \left( \cap _{n\geq 0}
g^{-q.n}(V_{i_0})\right) \cap (V_{i_1})^c\right] 
\end{equation}
is a non-empty, compact $g^q$-invariant subset of $D,$ contained in 
$$
TransDomain=\cup _{n\geq 0}g^{q.n}(V_{i_0}),\text{ an open set homeomorphic
to the plane.} 
$$

But $g^q$ is fixed point free in $V_{i_0}$ and $TransDomain$ is $g^q$%
-invariant, so 
$$
g^q\mid _{TransDomain}\text{ is a Brouwer homeomorphism.} 
$$
And this is a contradiction with the existence of the compact $g^q$%
-invariant set in (\ref{invset}). So the second case does not happen and
periodic prime ends are either accessible saddles or sources. $\Box $

\vskip 0.2truecm


\vskip 0.2truecm

\vskip 0.2truecm

\section{On $C^r$-generic families of diffeomorphisms}

Here we present results that characterize certain behaviors of 
$C^r$-generic one-parameter families of diffeomorphisms of surfaces. In
particular, one result is analogous to theorem \ref{comeconovo}, but now it
holds for residual subsets in the parameter space of $C^r$-generic
one-parameter families of diffeomorphisms.

\begin{theorem}
\label{generic1}\ (Brunovski) For any $r\in \{2,...,\infty \},$ if 
$(f_t)_{t\in I}$ is a $C^r$-generic one-parameter family of diffeomorphisms
of a closed Riemannian manifold, then periodic points are born from only two
different types of bifurcations: saddle-nodes and period doubling. In case
of saddle-nodes, if the parameter changes, then the saddle-node unfolds into a
saddle and a sink or source in one direction and the periodic point
disappears in the other direction. Moreover, at each fixed parameter, only
one saddle-node can happen.
\end{theorem}

{\bf Remark: }This statement appears in \cite{brunovski}. The exception is
the last part, which follows from the fact that the local Banach submanifold
in the space of $C^r$-diffeomorphisms, of maps which have a saddle-node
periodic orbit has codimension one. Therefore, a generic one-parameter
family of diffeomorphisms does pass through the intersection of two or more
such Banach manifolds. More parameters are necessary for that. See \cite
{sotom}.

\begin{theorem}
\label{generic2} For any $r\in \{1,2,...,\infty \},$ if $(f_t)_{t\in I}$ is
a $C^r$-generic one-parameter family of diffeomorphisms of a closed
orientable surface, then there exists a residual subset $R\subset I,$ such
that if $t\in R,$ then $f_t$ is a Kupka-Smale diffeomorphism.
\end{theorem}

{\bf Remark: }This theorem is folklore, but can be proved in the same way as
the next one, which is a version of theorem \ref{comeconovo} for families.

\begin{theorem}
\label{pixfami} For any $r\in \{1,2,...,\infty \},$ given a $C^r$-generic 
one-parameter family $(f_t)_{t\in I}$ of diffeomorphisms of any closed 
orientable surface,
 there exists a residual subset $R^{\prime }\subset
R\subset I$ $(R$ is the residual subset given by the previous theorem), such
that for $t\in R^{\prime },$ if $z_t$ is a hyperbolic $f_t$-periodic saddle
point with a full mesh, in case $w_t\in \overline{W^u(z_t)}$$,$ where $w_t$ is
another hyperbolic periodic saddle point for $f_t,$ then $W^u(z_t)$ actually
has a $C^1$-transverse intersection with $W^s(w_t).$ 
\end{theorem}

%
%
%
%
%

The proofs of theorems \ref{generic2} and \ref{pixfami} follow from the
Kupka-Smale theorem, theorem \ref{comeconovo} and the remark after it,
and the following general result:

\begin{lemma}
\label{andresva} If $r\geq 1$ and $R(M)$ is a $C^r$-residual subset of $%
\diff^r(M)$ (where $M$ is a compact manifold without boundary), then for a $%
C^r$-generic one-parameter family of diffeomorphisms $(f_t)_{t\in I}$, the
set of parameters $t\in I$ for which $f_t\in R(M)$ is residual in $I$.
\end{lemma}

{\it Proof:} Since $R(M)$ is residual, it contains some set of the form 
$$
R^{*}(M)=\bigcap_{i=1}^\infty W_n, 
$$
where each $W_n$ is open and dense in $\diff^r(M)$.

\begin{description}
\item[Claim: ]  For any fixed $s\in I$ and $m\in \mathbb{N}$, 
the set of one-parameter families $(f_t)_{t\in I}$ such that $f_s\in W_m$
is $C^r$-open and dense.
\end{description}

{\it Proof: }The openness is trivial. To prove the density, first recall
that the space of $C^r$-diffeomorphisms of a compact manifold $M$ is a
Banach manifold modeled on the space of $C^r$ vector fields $\mathfrak{X}%
^r(M)$. In fact, if $\exp {}_x:T_xM\to M$ denotes the exponential map for
any (smooth) complete Riemannian metric on $M$, a chart on a small enough
neighborhood $\mathcal{U}$ of $f$ is given by the map $\phi :\mathcal{U}%
\mapsto \mathfrak{X}^r(M)$ where $\phi (g)$ is the map $x\mapsto
exp_x^{-1}(g(f^{-1}(x)))$. Assuming $f=f_s$, there exists $\delta >0$ such
that $f_t\in \mathcal{U}$ whenever $|s-t|<\delta $. Let 
$\eta :\mathbb{R}\to \mathbb{R}$ be a $C^\infty $
function such that $\eta (t)=0$ if $|s-t|\geq \delta $ and $\eta (s)=1$. We
may choose a sequence $(g_n)_{n\in \mathbb{N}}$ of
elements of $W_m$ such that $g_n\to f$ in the $C^r$ topology. We further
assume that $\phi (\mathcal{U})$ is convex (by reducing it if necessary), so
we may define new maps by convex combination 
$$
f_{n,t}=\phi ^{-1}((1-\eta (t))\cdot \phi (f_t)+\eta (t)\cdot \phi (g_n)) 
$$
when $|s-t|<\delta $ and by $f_{n,t}=f_t$ otherwise. Note that $(x,t)\mapsto
f_{n,t}(x)$ is of class $C^r$ and it is straightforward to verify that $%
(f_{n,t})_{t\in I}\to (f_t)_{t\in I}$ in the $C^r$-topology as $n\to \infty $%
. 
This proves the claim. $\Box $

\vskip 0.2truecm

The lemma follows easily from the claim: letting $\{s_n:n\in \mathbb{N}\}$ be 
a dense subset of $I$ and applying the claim to $%
s=s_n$ and a given $m$ we obtain a $C^r$-open and dense set $\Lambda _{m,n}$
of $1$-parameter families $(f_t)_{t\in I}$ such that $f_{s_n}\in W_m$, and
thus 
$$
\Lambda =\bigcap_{n,m\in \mathbb{N}}\Lambda _{m,n} 
$$
is a residual set of families $(f_t)_{t\in I}$ such that $f_{s_n}\in W_m$
for all $n,m\in \mathbb{N}$. Note that this implies that
the set $W_m^{\prime }:=\{t\in I:f_t\in W_m\}$ is open and dense in $I,$
thus $R^{\prime }=\bigcap_{m\in \mathbb{N}}W_m^{\prime }$
is a residual subset of $I$ such that every $t\in R^{\prime }$ satisfies $%
f_t\in R^{*}(M)\subset R(M),$ completing the proof of the lemma. $\Box $

\vskip 0.2truecm

\section{Proofs of the main results}

In the first subsection, we prove an auxiliary result. 

\subsection{A general result}

\begin{proposition}
\label{geral} Assume $f\in \diff_0^{1+\epsilon }({\rm T^2})$ and its rotation
set has interior or $f\in \diff_0^{1+\epsilon }(S),$ where $S$ is a closed 
orientable surface of genus larger than $1,$ and it has a fully
essential system of curves {\bf $\mathscr{C}.$} If $D$ is a $f$-periodic
open disk, then either $\CDiam(D)<Max(f)$ (see proposition \ref{bounddiam}%
) or $D$ is homotopically unbounded. In case $D$ is homotopically unbounded, 
using the notation from
theorems \ref{casogeral}, \ref{casogeral2} and proposition \ref{bounddiam},
either $D$ avoids $\overline{W^s(Q)}$ and intersects $W^u(Q)$ (in this case $%
\partial D\supseteq \overline{W^s(Q)})$ or $D$ avoids $\overline{W^u(Q)}$
and intersects $W^s(Q)$ (in this case $\partial D\supseteq \overline{W^u(Q)}%
).$ Moreover, $\left( \lambda _u\cup \lambda _s\right) ^c\cap D=\cup
_{i=1}^\infty A_i,$ a disjoint union of open bounded disks, properly
labeled. If $n>0$ is a common period of $D,$ the point $Q$ and its stable and unstable branches,
and moreover, $\partial D\supseteq \overline{%
W^s(Q)},$ then $f^n(A_1)\subset A_1.$ All other $A_{i^{\prime }s}$ are
wandering and for each integer $i>1,$ there exists $m(i)$ such that $%
f^{n.m(i)}(A_i)\subset A_1.$ If $\partial D\supseteq \overline{W^u(Q)},$
then $f^n(A_1)\supset A_1.$ All other $A_{i^{\prime }s}$ are wandering and
for each integer $i>1,$ there exists $m(i)$ such that $f^{n.m(i)}(A_1)%
\supset A_i.$
\end{proposition}

{\it Proof:}

Without loss of generality, assume $D$ is $f$-invariant and $Q$ from
proposition \ref{bounddiam} is fixed, as are all four branches at $Q.$ If not,
consider some power of $f.$ Also, for simplicity of writing, we consider the
torus case.

If we lift $\lambda _u\cup \lambda _s$ to the plane (see proposition \ref
{bounddiam}), we get that its complement is made of disks whose diameters
are smaller than $Max(f).$ Therefore if some connected component $\widetilde{%
D}$ of $\pi $$^{-1}(D)$ intersects $\pi $$^{-1}(\lambda _u\cup \lambda _s),$
then $\widetilde{D}$ contains an arc $\alpha $ transversal to $W^u( 
\widetilde{Q})$ or $W^s(\widetilde{Q}).$ Choose some $(p,r)$ such that $
\widetilde{f}(\widetilde{D})-(p,r)=\widetilde{D}.$ The $\lambda $-lemma 
then implies that the orbit of this arc under $\widetilde{f}(\bullet
)-(p,r)$ is unbounded and so $\widetilde{D}$ is unbounded.

If $D$ intersects $W^u(Q)$ and $W^s(Q),$ as it is open, $D$ contain arcs $%
\gamma _u$ and $\gamma _s,$ $C^1$-transversal, respectively, to $W^s(Q)$ and 
$W^u(Q).$ So, as $n>0$ goes to infinity, $f^n(\gamma _u)$ $C^1$-accumulates
on $W^u(Q),$ the same for $f^{-n}(\gamma _s)$ with respect to $W^s(Q).$ As 
$W^u(Q)\cap W^s(Q)$ contains a topologically transverse point $z,$ for which 
the arc in $W^u(Q)$
whose endpoints are $Q$ and $z,$ union with the arc in $W^s(Q)$ whose
endpoints are $Q$ and $z$ form a homotopically non-trivial closed curve, 
$f^n(\gamma _u)\cup f^{-n}(\gamma _s) $ also contains a homotopically
non-trivial closed curve for all sufficiently large $n>0$ (this follows from
the topological transversality between $W^s(Q)$ and $W^u(Q)$ at $z$). But as 
$f^n(\gamma _u)\cup f^{-n}(\gamma _s)\subset D$ for all $n>0,$ this is a
contradiction with the fact that $D$ is a disk.

Clearly, if $D$ is unbounded and avoids $\overline{W^s(Q)},$ then it
intersects $W^u(Q).$ So, it contains an arc whose negative iterates $C^1$%
-accumulate on $W^s(Q).$ Therefore, $\partial D\supseteq \overline{W^s(Q)}.$

Finally, let us look at the connected components of $\left( \lambda _u\cup
\lambda _s\right) ^c\cap D,$ assuming that $\partial D\supseteq \overline{%
W^s(Q)}.$ Each one is an open disk, whose diameter is bounded by $Max(f)>0$
(see proposition \ref{bounddiam}). Choose any point $z^{*}\in D.$ Let $\beta $
be an arc in $D$ whose endpoints are $z^{*}$ and $f(z^{*}).$ As $%
f^{-1}(\lambda _u)\subset \lambda _u$ and $Q\notin D,$ there exists $m(\beta
)>0$ such that for all $i\geq m(\beta ),$ $f^i(\beta )$ avoids $\lambda
_u\cup \lambda _s.$ Therefore, for all $i\geq m(\beta ),$ $f^i(z^{*})$
belongs to the same connected component of $\left( \lambda _u\cup \lambda
_s\right) ^c\cap D.$ Call it $A_1.$

Consider any other point $z$ in $D.$ Fix an arc $\alpha \subset D,$
connecting $z^{*}$ to $z.$ As for $\beta ,$ there exists $m(\alpha )>m(\beta
)>0$ such that for all $i\geq m(\alpha ),$ $f^i(\alpha )$ avoids $\lambda
_u\cup \lambda _s.$ Therefore, for all $i\geq m(\alpha ),$ $f^i(\alpha )$ is
contained in $A_1.$ So the positive orbit of any $z\in D,$ at some point
enters $A_1$ and stays there forever.

Moreover, if an arc $\gamma \subset D$ avoids $\lambda _u,$ then $f^i(\gamma
)$ avoids $\lambda _u$ for all $i>0.$ This implies that the image under $f$
of any connected component of $\left( \lambda _u\cup \lambda _s\right)
^c\cap D$ is always contained in another connected component. So
immediately, $f(A_1)\subset A_1.$ And for the others? Fix some $i>1$ and
look at $A_i.$ If for some $m>0,$ $f^m(A_i)\cap A_i\neq \emptyset ,$ then $%
f^m(A_i)\subset A_i\Rightarrow f^{j.m}(A_i)\subset A_i$ for all integers $%
j>0.$ And this is a contradiction with the fact that all points in $D$ enter 
$A_1$ after a certain positive iterate and do not leave it anymore.

So for all $i>1,$ each $A_i$ is wandering and for a certain $m(i)>0,$ $%
f^{m(i)}(A_i)$ is contained in $A_1.$

The case when D is unbounded and intersects $W^s(Q)$ is analogous. $\Box $

\vskip 0.2truecm

\subsection{Proof of theorem 1} 

The proof of this theorem follows entirely from theorem \ref{repulsor0}. 
If $D$ is a 
$f$-periodic open disk which is a connected component of $\left( \overline{%
W^s(z)}\right) ^c$ for some saddle $z$ with a full mesh, and the prime ends
rotation number of $D$ is rational, then theorem \ref{repulsor0} implies that 
$\widehat{D}$ is a closed disk in
whose boundary there are only finitely many periodic prime ends.
Necessarily, some are accessible hyperbolic saddles and the remaining are
sources. For such an accessible saddle $w\in \partial D,$ the whole stable
manifold of $w$ is contained in $\partial D,$ some unstable branch at 
$w$ has a transversal intersection with $W^s(z)$ and one unstable branch at $w$
is contained in $D.$ So, $W^s(w)$ is accessible from $D.$ 
Clearly, sources in $\partial \widehat{D}$ might correspond to accessible
hyperbolic periodic sources in $\partial D,$ but, necessarily, 
in case $D$ is homotopically unbounded, in any given neighborhood of 
some of these sources is "hidden the unboundedness of $D".$ By this, we mean that
for some sources in $\partial \widehat{D},$ and any crosscut around one of them, the 
corresponding cross-section in $D$ is an homotopically unbounded disk. Moreover, 
$D$ minus the union 
of these particular homotopically unbounded
 cross-sections is always a homotopically bounded disk. $\Box $ 
\vskip 0.2truecm

\subsection{Proof of theorem C}

This is immediate from the dynamics of $\widehat{f}^n\mid _{\partial 
\widehat{D}},$ see figure 3. $\Box$

\vskip0.2truecm

\begin{figure}[!h]
	\centering
	\includegraphics[scale=0.53]{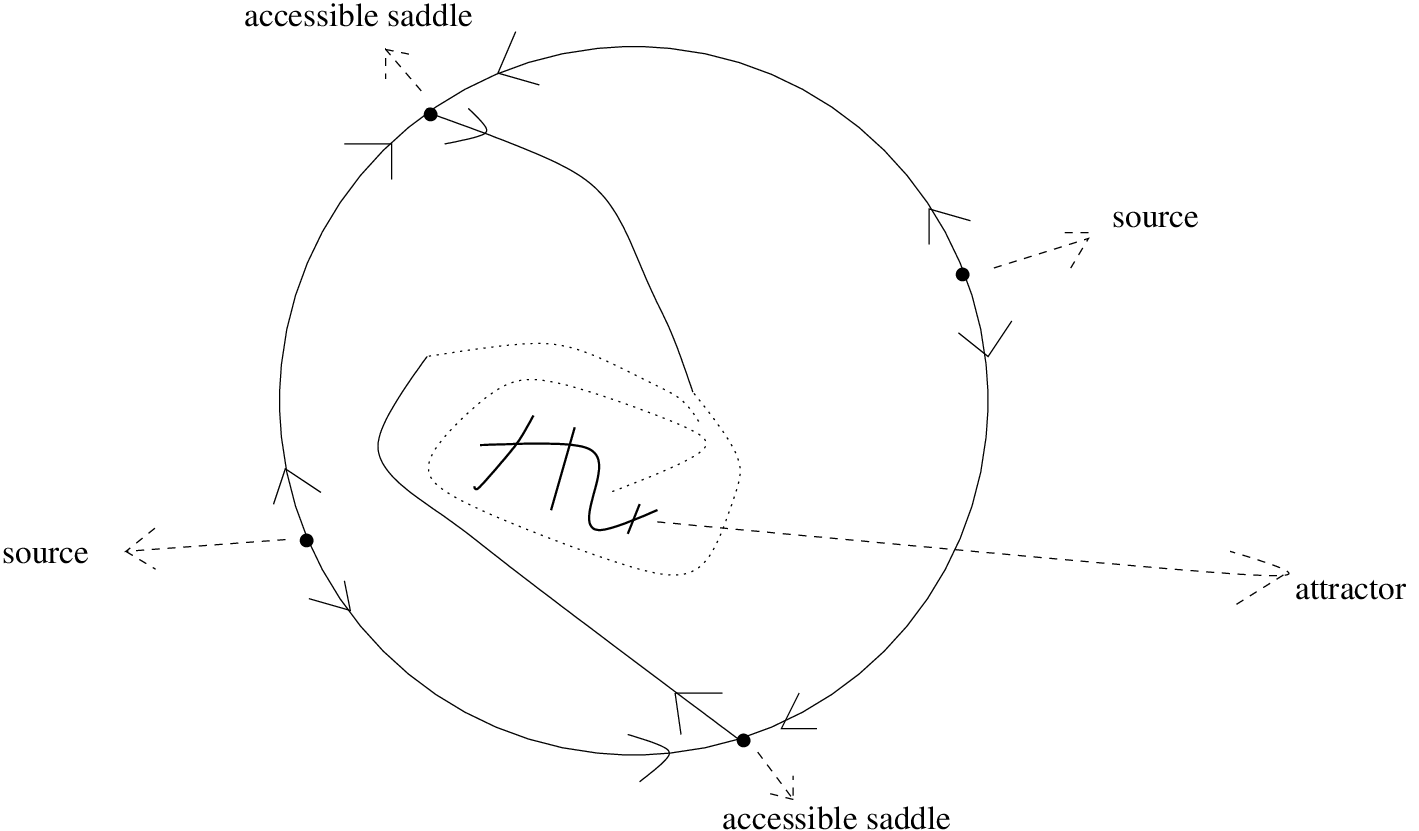}
	\caption{Diagram showing the dynamics of $\widehat{f}^n$ restricted to $\widehat{D}.$}
	\label{clopen4}
\end{figure}

\vskip 0.2truecm

\subsection{Proof of theorems B and B'} 

Any homotopically unbounded, maximally periodic open disk $D$ must be a connected component 
of the complement of $\overline{W^s(z)}$  (or of $\overline{W^u(z)}$), where $z$ is a contractible hyperbolic periodic 
saddle with a full mesh. This follows from the following: as $D$ is unbounded, it intersects
stable or unstable manifolds of $z.$ From proposition \ref{geral}, $D$ can not intersect both. 
$\Box$

\vskip0.2truecm

\subsection{Proof of theorem D}

Without loss of generality, 
assume $D$ is a $f^N$-periodic (for some integer 
$N>0$ which is a multiple of $2.n$) connected component of 
$\left( \cup _{i=0}^{n-1}\overline{W^s(f^i(z))}\right) ^c.$ 
The choice of $N$ implies that $f^N$ leaves each boundary component of $D$ invariant. 
Also, the prime ends rotation number of $f^N$ restricted to each component of $\partial D$ 
is assumed to be rational. From observation 2 before the statement of theorem 
\ref{repulsor0}, the prime ends rotation number of $f^N$ 
is the same at all components of $\partial D$. 

If instead
of stable manifold at the orbit of $z,$ it were unstable manifold, we just
would have to work with $f^{-1}.$

From part B) of theorem \ref{repulsor0},
each connected component of the boundary of $\widehat{D},$ the prime ends
compactification of $D,$ is a circle with only finitely many periodic prime
ends. Necessarily, some are accessible hyperbolic saddles and the remaining
are sources. For such an accessible saddle $w\in \partial D,$ 
$\overline{W^s(w)}$ is contained in $\partial D$ and one 
unstable branch, denoted $\lambda _w^u,$ is contained in $D.$
%
%
So, a result analogous to theorem 1 holds for each connected component $\Gamma$ of $\partial D.$ 
The rest of the
proof follows from the dynamics of $\widehat{f}^N\mid _{\widehat{\Gamma}},$ again 
see figure 3. $\Box $
\vskip0.2truecm

\subsection{Proof of theorems A and A'}
   
From the theorem hypotheses, there exists $t_0<t^{*}<t_1$ such that $%
f_{t^{*}}$ only has one $q$-periodic orbit (for some integer $q>0)$ of
rotation vector $\rho $$,$ which is a saddle-node (see theorem \ref{generic1}%
). And for some $\epsilon >0,$ for all $t\in ]t^{*},t^{*}+\epsilon [,$ $f_t$
has exactly two $q$-periodic orbits of rotation vector $\rho $$,$ a saddle
and a sink (or a source), one branch of the saddle connected to the sink (or
source). Moreover, if $t<t^{*},$ then $f_t$ has no $q$-periodic points of
rotation vector $\rho,$ see \cite{palistakens} page 38-39 and \cite{carr}. This 
is a point where we need $C^2$ differentiability.

Assume without loss of generality that $q=1$ and the saddle-node unfolds for 
$t>t^{*}$ into a saddle $p_t$ and a sink $s_t.$ Let $Fix_\rho (f_t)=\{z\in
Fix(f_t):$ $z$ has rotation vector $\rho \}$ and let $\alpha _t^u$ be the
unstable branch of $p_t$ that connects it to $s_t.$ If $D_t$ is the basin of
attraction of $s_t,$ we are going to show that for $t$ in an open subset of $%
]t^{*},t^{*}+\epsilon [,$ $D_t$ is an homotopically unbounded $f_t$-invariant disk.
Clearly, the prime ends compactification of $D_t$ has zero rotation number,
because $p_t$ is a saddle in $\partial D_t,$ accessible through $\alpha
_t^u, $ which is $f_t$-invariant. This means that $W^s(p_t)\subset \partial
D_t.$

Either when the surface is the torus or has higher genus, our assumptions imply that
for all $t\in I,$ $f_t$ has contractible 
periodic saddles with a full mesh. So in both cases, $f_t$ has 
a periodic saddle $R_t$ with a
full mesh for all $t\in [t^{*},t^{*}+\epsilon ],$ (this is clearly true if $%
\epsilon >0$ is sufficiently small because transverse intersections persist
under perturbations and the existence of certain finitely many transverse
intersections imply the existence of a full mesh). The complement of $
\overline{W^u(R_t)}$ is a union of open topological disks. From theorem \ref
{pixfami}, there exists a residual subset $R^{\prime }\subset
]t^{*},t^{*}+\epsilon [$ such that for any $t\in R^{\prime }:$

\begin{enumerate}
\item  If $p_t\in \overline{W^u(R_t)},$ then $W^u(R_t)$ has a $C^1$%
-transverse intersection with $W^s(p_t).$ But then, $\overline{W^s(R_t)}%
\subset \overline{W^s(p_t)}$ and so $\partial D_t\supset \overline{W^s(p_t)}%
\supset \overline{W^s(R_t)},$ which implies that $D_t$ is unbounded.

\item  If $p_t\notin \overline{W^u(R_t)},$ then there exists a $f_t$%
-invariant connected component $M_t$ of $\left( \overline{W^u(R_t)}\right)
^c $ which contains $p_t.$ As we said before, it is an open disk. As 
$p_t\in M_t$, it is easy to see that $\alpha _t^u$ is also 
contained in $M_t.$ There are two cases:

\begin{itemize}
\item  $s_t\in \partial M_t.$ This implies that $D_t$ intersects $W^u(R_t),$
and so, negative iterates of $D_t$ converge to $W^s(R_t),$ thus $D_t$ is unbounded.

\item  $\{p_t\cup \alpha _t^u\cup s_t\}\subset M_t.$ When this happens,
deform $f_t$ inside a small neighborhood $V$ of $p_t\cup \alpha _t^u\cup
s_t, $ $closure(V)\subset M_t$ in order to get a diffeomorphism $g$ which
coincides with $f_t$ outside $V,$ with the property that in the whole
surface, $g$ has no fixed points of rotation vector $\rho .$ As $closure(V)$
avoids $\overline{W^u(R_t)},$ it also avoids a local stable manifold at $R_t$
for $f_t.$ So, although the stable manifold of $R_t$ for the map $g$ may be
different from the stable manifold of $R_t$ for $f_t,$ $R_t$ still has a
full mesh for $g,$ $g(M_t)=M_t$ and $\partial M_t\subseteq \overline{W^u(R_t)%
}.$ As $M_t$ is a $g$-invariant open disk of rotation vector $\rho ,$ if it
is bounded, then $\overline{M_t}$ contains a fixed point of rotation vector $%
\rho .$ If $M_t$ is unbounded, proposition \ref{geral} implies the existence
of a bounded closed disk $\overline{A_1}\subset \overline{M_t}$ which
satisfies $g^{-1}(\overline{A_1})\subseteq \overline{A_1},$ so $\overline{A_1%
}$ must also have fixed points (of rotation vector $\rho $). This
contradiction shows that it is not possible that $\{p_t\cup \alpha _t^u\cup
s_t\}\subset M_t.$ %
\end{itemize}
\end{enumerate}

So, either $D_t$ intersects $W^u(R_t),$ which is an open condition; or $%
\partial D_t\supset \overline{W^s(p_t)}\supset \overline{W^s(R_t)}.$ As the
last relation comes from the fact that $W^s(p_t)$ has $C^1$-transverse intersections
with $W^u(R_t),$ also an open condition, the proof is over. $\Box $

\end{document}